\documentclass[11pt]{amsart}
\usepackage{amsmath}
\usepackage{amssymb}
\usepackage{amsfonts}
\usepackage{amsthm}
\usepackage{enumerate}
\usepackage[mathscr]{eucal}
\usepackage{verbatim}
\usepackage{amscd}
\usepackage{appendix}
\usepackage{tikz}
\usepackage{hyperref}
\usepackage{cleveref}

\numberwithin{equation}{section}

\newtheorem{innercustomgeneric}{\customgenericname}
\providecommand{\customgenericname}{}

\newcommand{\newcustomtheorem}[2]{\newenvironment{#1}[1]
  {\renewcommand\customgenericname{#2}
   \renewcommand\theinnercustomgeneric{##1}\innercustomgeneric}{\endinnercustomgeneric}}

\newcustomtheorem{customthm}{Theorem}

\newcommand{\newcustomlemma}[2]{\newenvironment{#1}[1]
  {\renewcommand\customgenericname{#2}
   \renewcommand\theinnercustomgeneric{##1} \innercustomgeneric}{\endinnercustomgeneric}}

\newcustomlemma{customlemma}{Lemma}

\newcustomlemma{customproposition}{Proposition}

\newcustomlemma{customclaim}{Claim}

\newcustomlemma{customdef}{Definition}

\theoremstyle{plain}
\newtheorem{theorem}{Theorem}
\newtheorem{corollary}[theorem]{Corollary}
\newtheorem{lemma}[theorem]{Lemma}
\newtheorem{proposition}[theorem]{Proposition}

\newtheorem*{theorem*}{Theorem}
\newtheorem*{lemma*}{Lemma}
\newtheorem*{proposition*}{Proposition}
\newtheorem*{corollary*}{Corollary}
\newtheorem*{remark*}{Remark} 
\newtheorem*{remarks*}{Remarks}
\newtheorem*{conj*}{Conjecture}

\newcommand{\LL}{\mathcal{L}}

\newcommand{\bbz}{\mathbb{Z}}

\newcommand{\bbrn}{\mathbb R^n}

\newcommand{\bbc}{\mathbb{C}}

\newcommand{\xxxi}{\vec{\boldsymbol{\xi}\;}}

\newcommand{\qqq}{\vec{\boldsymbol{q}}}
\newcommand{\www}{\vec{\boldsymbol{w}}}
\newcommand{\yyy}{\vec{\boldsymbol{y}}}
\newcommand{\zzz}{\vec{\boldsymbol{z}}}
\newcommand{\uuu}{\vec{\boldsymbol{u}}}

\newcommand{\ppp}{\vec{\boldsymbol{p}}}

\def\000{\vec{\boldsymbol{0}}}

\def\ga{\gamma}

\def\la{\lambda}

\def\Om{\Omega}

\newcommand{\q}{\quad}
\newcommand{\qq}{\qquad}

\DeclareFontFamily{U}{mathx}{\hyphenchar\font45}
\DeclareFontShape{U}{mathx}{m}{n}{
	<5> <6> <7> <8> <9> <10>
	<10.95> <12> <14.4> <17.28> <20.74> <24.88>
	mathx10
}{}

\def\wh{\widehat}

\newcommand{\wt}{\widetilde}

\makeatletter
\@namedef{subjclassname@2020}{\textup{2020} Mathematics Subject Classification}
\makeatother

\allowdisplaybreaks

\makeatletter
\@namedef{subjclassname@2020}{\textup{2020} Mathematics Subject Classification}
\makeatother

\makeindex         

\begin{document}

\author{Bae Jun Park}
\address{B. Park, Department of Mathematics, Sungkyunkwan University, Suwon 16419, Republic of Korea}
\email{bpark43@skku.edu}

\thanks{The author is supported in part by NRF grant RS-2025-20512969 and by POSCO Science Fellowship of POSCO TJ Park Foundation. The author is grateful for support by the Open KIAS Center at Korea Institute for Advanced Study.}

 \title{Weighted Estimates for Multilinear Singular Integrals with Rough Kernels} 

\begin{abstract} 
We establish weighted norm inequalities for multilinear singular integral operators with rough kernels. Specifically, we consider the multilinear singular integral operator $\mathcal{L}_\Omega$ associated with an integrable function $\Omega$ on the unit sphere $\mathbb{S}^{mn-1}$ satisfying the vanishing mean condition. Extending the classical results of Watson \cite{Wa1990} and Duoandikoetxea \cite{Duan1993} to the multilinear setting, we prove that $\mathcal{L}_\Omega$ is bounded from $L^{p_1}(w_1)\times\cdots\times L^{p_m}(w_m)$ to $L^p(v_{\www})$ under the assumption that $\Omega\in L^q(\mathbb{S}^{mn-1})$ and that the $m$-tuple of weights $\www= (w_1,\ldots,w_m)$ lies in the multiple weight class $\mathrm{A}_{\ppp/q'}$. Here, $q'$ denotes the H\"older conjugate of $q$, and we assume $q'\le p_1,\dots,p_m<\infty$ with $1/p = 1/p_1 + \cdots + 1/p_m$.
\end{abstract}

\maketitle


\section{Introduction}

Let $\Omega$ be an integrable function on the unit sphere $\mathbb{S}^{n-1}$ with mean value zero. The associated singular integral operator $T_{\Omega}$ is then defined by
$$T_{\Omega}f(x):=\mathrm{p.v.}\int_{\bbrn}K_{\Omega}(y)f(x-y)\; dy$$
where $K_{\Omega}(y):=\frac{\Omega(y/|y|)}{|y|^n}$.
This class of operators was introduced by Calder\'on and Zygmund \cite{Ca_Zy1952} and  has been a central object of study in harmonic analysis, attracting considerable attention,  since the pioneering work of Calder\'on and Zygmund \cite{Ca_Zy1956} who proved the $L^p$ boundedness for $T_{\Omega}$ under the assumption $\Omega\in L\log{L}(\mathbb{S}^{n-1})$. The condition $\Omega\in L\log{L}(\mathbb{S}^{n-1})$ was later relaxed by Coifman and Weiss \cite{Co_We1977} and Connett \cite{Co1979} who established the $L^p$ boundedness for $T_{\Omega}$ if $\Omega$ belongs to the Hardy space $H^1(\mathbb{S}^{n-1})$. As an endpoint estimate, Christ and Rubio de Francia \cite{Ch_Ru1988} extended the previous $L^p$ boundedness results to the weak type $(1,1)$ boundedness when $n=2$ and $\Omega\in L\log{L}(\mathbb{S}^1)$, and independently the same conclusion was also obtained by Hofmann \cite{Ho1988}. Finally, Seeger \cite{Se1996} extended the weak-type $(1,1)$ boundedness to all dimensions $n\ge 2$ assuming $\Omega\in L\log L(\mathbb{S}^{n-1})$.
We also refer to \cite{St2001, Ta1999} for further results.
 
 Alongside these unweighted results, considerable progress has been made in understanding weighted norm inequalities for rough singular integrals.
  Recall from \cite{Mu1972} that a nonnegative, locally integrable function $w$ on $\bbrn$ is said to belong to the Muckenhoupt $A_p$ class for $1\le p<\infty$ if the associated $A_p$ constant, denoted by $[w]_{A_p}$, is finite
 where 
 $$[w]_{A_1}:= \sup_{Q:\text{cubes in $\bbrn$}}\Big(\frac{1}{|Q|}\int_Q w(x)\; dx\Big)\Big( \inf_{x\in Q}{w(x)}\Big)^{-1},$$
\begin{align*}
[w]_{A_p}:= \sup_{Q:\text{cubes in $\bbrn$}}\Big( \frac{1}{|Q|}\int_Q w(x)\; dx\Big)\Big(\frac{1}{|Q|}\int_Q \big( w(x)\big)^{-\frac{1}{p-1}}\; dx \Big)^{p-1},\q 1<p<\infty.
 \end{align*}
Here, $\mathcal{M}$ is the Hardy-Littlewood maximal operator, whose definition will be recalled in Section \ref{maxinesec}. 
 For $p=\infty$, we define $A_{\infty}:=\bigcup_{p>1}A_p$ and then it turns out that 
 \begin{equation}\label{apweightinclusion}
 A_p\subset A_q \q  \text{ if}~ 1\le p\le q\le \infty.
 \end{equation}
Given a weight $w$,
the weighted Lebesgue space $L^p(w)$, $0<p<\infty$, consists of all measurable functions $f$ on $\bbrn$ satisfying
$$\|f\|_{L^p(w)}:=\left(\int_{\bbrn}|f(x)|^p w(x)\, dx \right)^{1/p}<\infty . $$
 Based on the foundational theory of Muckenhoupt $A_p$ weights, Duoandikoetxea and Rubio de Francia \cite{Du_Ru1986} showed that if $w\in A_p$, then 
  $T_{\Omega}$ is bounded on $L^p(w)$ when  $\Omega\in {L^{\infty}(\mathbb{S}^{n-1})}$. 
  These estimates were later refined by Watson \cite{Wa1990} and Duoandikoetxea \cite{Duan1993}.
 Let $\mathscr{C}_c^{\infty}(\bbrn)$ denote the family of all smooth functions with compact support on $\bbrn$.
\begin{customthm}{A}\cite{Duan1993, Wa1990}\label{weightlinear}
Let $\Omega$ be an integrable function on $\mathbb{S}^{n-1}$ with mean value zero.
Let $1<p<\infty$ and $1<q\le \infty$ satisfy $q'\le p$.
Suppose that $\Omega\in L^q(\mathbb{S}^{n-1})$ and $w\in A_{p/q'}$.
Then there exists a constant $C$, depending on the $A_{p/q'}$-constant $[w]_{A_{p/q'}}$, such that
\begin{equation*}
\big\Vert T_{\Omega}f\big\Vert_{L^p(w)}\le C \Vert \Omega\Vert_{L^q(\mathbb{S}^{n-1})}\Vert f\Vert_{L^p(w)}
\end{equation*}
for all $f\in\mathscr{C}_c^{\infty}(\bbrn)$.
\end{customthm}
 We also refer to \cite{Va1996} for weighted weak type $(1,1)$ estimates, and \cite{Co_Cu_Di_Ou2017, Hy_Ro_Ta2017, Li_Pe_Ri_Ro2019} for quantitative weighted bounds for the operator $T_{\Omega}$.

\hfill

The purpose of this paper is to provide a multilinear extension of Theorem \ref{weightlinear} with multiple weights, introduced by Lerner, Ombrosi, P\'erez, Torres, and Trujillo-Gonz\'alez \cite{Le_Om_Pe_To_Tr2009}.
Let $m$ be a positive integer greater than $1$, and let $\Omega$ be now an integrable function on the unit sphere $\mathbb{S}^{mn-1}$ having mean value zero 
\begin{equation*}
\int_{\mathbb{S}^{mn-1}}\Omega(\yyy')~ d\sigma(\yyy')=0
\end{equation*} 
where $d\sigma$ stands for the surface measure on $\mathbb{S}^{mn-1}$, $\yyy=(y_1,\dots,y_m)\in (\bbrn)^m$,  and $\yyy':=\frac{\yyy}{|\yyy|}\in \mathbb{S}^{mn-1}$.
We define a kernel
\begin{equation*}
K(\yyy):=\frac{\Omega(\yyy')}{|\yyy|^{mn}}, \qquad \yyy \neq \000.
\end{equation*} 
Then the corresponding multilinear singular integral operator $\LL_{\Omega}$ is defined by
\begin{align*}
\mathcal L_{\Om}\big(f_1,\dots,f_m\big)(x)& :=\mathrm{p.v.} \int_{(\bbrn)^m}{K(\yyy) f_1(x-y_1)\cdots f_m(x-y_m)}~d\yyy
\end{align*}
 for  $f_1,\dots,f_m\in\mathscr{C}_c^{\infty}(\bbrn)$.

The boundedness properties of bilinear singular integral operators in the one-dimensional setting ($n=1$) were first investigated by Coifman and Meyer \cite{Co_Me1975}, who proved $L^{p_1}\times L^{p_2}\to L^p$ boundedness when $\Omega$ is of bounded variation on the unit circle $\mathbb{S}^1$. 
Grafakos and Torres \cite{Gr_To2002} extended this result to higher dimensions  and multilinear settings, assuming Lipschitz regularity of $\Omega$. Both of these results relied on certain smoothness assumptions on the angular part $\Omega$ of the kernel.
Subsequent progress in the bilinear case was made by Grafakos, He, and Honz\'ik \cite{Gr_He_Ho2018}, who considered the case when $\Omega$ is merely bounded, i.e., $\Omega \in L^\infty(\mathbb{S}^{2n-1})$. In particular, they established the initial estimates $L^2 \times L^2 \to L^1$ for $\LL_{\Omega}$, even when $\Omega$ belongs to $L^2(\mathbb{S}^{2n-1})$, by employing a novel approach based on Daubechies wavelet decomposition \cite{Da1988}. This result was further refined by Grafakos, He, and Slav\'ikov\'a \cite{Gr_He_Sl2020}, who weakened the assumption on $\Omega$ to $\Omega \in L^q(\mathbb{S}^{2n-1})$ for $q > \frac{4}{3}$. A more general boundedness result was obtained by He and the author \cite{He_Park2023}, who extended  the range of exponents  to all $1 < p_1, p_2 < \infty$ and $\frac{1}{2} < p < \infty$ with $1/p=1/p_1+1/p_2$, under the assumption that $\Omega \in L^q(\mathbb{S}^{2n-1})$ for $q > \max\left(\frac{4}{3}, \frac{p}{2p-1}\right)$.
In the general multilinear setting, Grafakos, He, Honz\'ik, and the author \cite{Gr_He_Ho_Park2023} obtained an initial boundedness result for $\LL_{\Omega}$ from $L^2 \times \cdots \times L^2$ to $L^{2/m}$ when $\Omega \in L^q(\mathbb{S}^{mn-1})$ for $q > \frac{2m}{m+1}$. The proof still relied heavily on the wavelet decomposition of Daubechies, although more delicate technical challenges arose due to the fact that the target space $L^{2/m}(\mathbb{R}^n)$ is not a Banach space when $m \ge 3$. This multilinear boundedness was later extended to the full range $1 < p_1, \dots, p_m < \infty$ and $\frac{1}{m} < p < \infty$ in \cite{Gr_He_Ho_Park_JLMS}. 
More recently, Dosidis and Slav\'ikov\'a \cite{Do_Sl2024} improved these results in certain regimes of exponents, and in particular showed that the minimal assumption $\Omega \in L^q(\mathbb{S}^{mn-1})$ with $q > 1$ is sufficient for the boundedness of $\LL_{\Omega}$ from $L^{p_1} \times \cdots \times L^{p_m}$ to $L^p$ whenever $1 < p, p_1, \dots, p_m < \infty$ and $1/p=1/p_1+\cdots+1/p_m$.
We present most of these results in the following one formulation:
\begin{customthm}{B}\cite{Do_Sl2024, Gr_He_Ho2018, Gr_He_Ho_Park2023, Gr_He_Sl2020, He_Park2023}\label{knownbdresult}
Let $1<p_1,\dots,p_m<\infty$ and $\frac{1}{m}<p<\infty$ with ${1}/{p}={1}/{p_1}+\cdots+{1}/{p_m}$.
Suppose that $1<q\le \infty$ satisfies
\begin{equation}\label{qcondequi}
\sum_{j\in J}\frac{1}{p_j}<\frac{|J|}{q'}+\frac{1}{q} \q \text{ for every subset $J$ of $\{1,\dots,m\}$}.
\end{equation}
 Then there exists a constant $C>0$ such that
\begin{equation*}
\big\Vert \LL_{\Omega}(f_1,\dots,f_m)\big\Vert_{L^p(\bbrn)}\le C\Vert \Omega\Vert_{L^{q}(\mathbb{S}^{mn-1})}\prod_{j=1}^{m}\Vert f_j\Vert_{L^{p_j}(\bbrn)}
\end{equation*}
for all $f_1,\dots,f_m\in\mathscr{C}_c^{\infty}(\bbrn)$.
\end{customthm}

\hfill

To formulate weighted estimates for multilinear operators, we employ the multiple weight theory developed by Lerner, Ombrosi, P\'erez, Torres, and Trujillo-Gonz\'alez \cite{Le_Om_Pe_To_Tr2009}.
 \begin{customdef}{C}\cite{Le_Om_Pe_To_Tr2009}
 Let $1\le p_1,\dots,p_m<\infty$ and $1/p=1/p_1+\dots+1/p_m$.
 Then we define $\mathrm{A}_{\ppp}$, $\ppp=(p_1,\dots,p_m)$, to be the class of $m$-tuples of weights $\www:=(w_1,\dots,w_m)$ satisfying
 \begin{equation*}
 \sup_{Q: \text{$\mathrm{cubes}$ $\mathrm{in}$ $\bbrn$}}\bigg[ \Big( \frac{1}{|Q|}\int_Q v_{\www}(x)\; dx \Big)^{{1}/{p}}\prod_{j=1}^{m}\Big(\frac{1}{|Q|}\int_Q \big(w_j(x)\big)^{1-p_j'}\; dx \Big)^{{1}/{p_j'}}\bigg]<\infty
 \end{equation*}
 where $p_j'$ denotes the H\"older conjugate of $p_j$ and
 $$v_{\www}(x):=\prod_{j=1}^{m}\big( w_j(x) \big)^{{p}/{p_j}}.$$
 When $p_j=1$, $(\frac{1}{|Q|}\int_Q w_j^{1-p_j'})^{1/p_j'}$ is understood as $(\inf_Q w_j)^{-1}$.
 \end{customdef}
 The multiple weight space serves as a natural extension of the classical Muckenhoupt weight class, as it coincides with $A_p$ when $m=1$, and it also has a maximal function characterization analogous to that of the $A_p$. See \eqref{chaweight} and Lemma \ref{multiweightcha} below.
 However, it does not necessarily preserve all the properties of the Muckenhoupt weights, which makes the theory of multiple weights more intricate.
 For example, as mentioned in \cite[Remark 7.3]{Le_Om_Pe_To_Tr2009}, the classes $\mathrm{A}_{\ppp}$ are not generally increasing with the natural partial order, unlike \eqref{apweightinclusion}. Further discussion on the properties of $\mathrm{A}_{\ppp}$ will be presented in the next section.

The main result of this paper is the following weighted estimate for multilinear rough singular integrals.
\begin{theorem}\label{singularthm}
Let $1<p_1,\dots,p_m<\infty$, and $\frac{1}{m}<p<\infty$ with $\frac{1}{p}=\frac{1}{p_1}+\cdots+\frac{1}{p_m}$.
Suppose that $1<q\le \infty$ satisfies $ q'\le p_1,\dots,p_m$ and $\www\in \mathrm{A}_{(p_1/q',\dots,\,p_m/q')}$.
Then we have
$$\big\Vert \LL_{\Omega}(f_1,\dots,f_m )\big\Vert_{L^p(v_{\www})}\lesssim \Vert \Omega\Vert_{L^q(\mathbb{S}^{mn-1})}\prod_{j=1}^{m}\Vert f_j\Vert_{L^{p_j}(w_j)}$$
for all $f_j \in\mathscr{C}_c^{\infty}(\bbrn)$.
\end{theorem}
We remark that if $ q'\le p_1,\dots,p_m$ and $1<p_1,\dots,p_m$, then \eqref{qcondequi} holds.\\

The proof of Theorem \ref{singularthm} is based on a dyadic decomposition of the kernel $K$ of  $\mathcal{L}_\Omega$ and kernel estimates established by Duoandikoetxea and Rubio de Francia \cite{Du_Ru1986}. We begin by decomposing the kernel $K$ into localized components and further decomposing each component via a Littlewood-Paley decomposition, which leads to a double-indexed kernel decomposition denoted by $K_\mu^\gamma$. This allows us to control the size and smoothness of each piece precisely. The resulting analysis naturally separates into two parts: the low-frequency part ($\mu \le 0$) and the high-frequency part ($\mu > 0$), due to their substantially different analytic behaviors.
For the low frequency part, the corresponding kernel satisfies the standard size and smoothness conditions of multilinear Calder\'on-Zygmund operators. Thus, the multilinear weighted theory developed in \cite{Le_Om_Pe_To_Tr2009} can be applied, yielding the desired weighted bounds under the assumption that the multiple weights belong to the class $\mathrm{A}_{\ppp}$. We note that the assumption $\vec{w} \in \mathrm{A}_{\ppp/q'}$ in Theorem \ref{singularthm} ensures this inclusion, as $\mathrm{A}_{\ppp/q'} \subset \mathrm{A}_{\ppp}$ (see Lemma \ref{mweightincl}).
In contrast, the high frequency part requires a different approach. Although each $K_\mu^\gamma$ again satisfies the multilinear Calder\'on-Zygmund kernel conditions, the constants in the associated size and smoothness bounds exhibit exponential growth in $\mu>0$. As a result, the multilinear Calder\'on-Zygmund theory cannot be applied directly, since the resulting estimates are not summable over $\mu > 0$.
To overcome this difficulty, we establish a pointwise estimate for the sharp maximal function associated with each high frequency piece (see Proposition \ref{mainpropo}). This allows us to derive a new weighted estimate with only polynomial growth in $\mu$ (see Corollary \ref{qconjulessp}). We then apply a multilinear version of Stein's complex interpolation theorem, which also accommodates interpolation of weights. Specifically, by interpolating between an unweighted multilinear estimate with exponential decay in $\mu>0$, available from previous results, and the new weighted estimate with polynomial growth, we obtain sufficient exponential decay in $\mu$ to ensure summability over $\mu>0$.

\hfill

The paper is organized as follows.
In Section \ref{prelisec}, we present several preliminary results including maximal inequalities, fundamental properties of multiple weights, and complex interpolation for analytic families of multilinear operators. Section \ref{prthm1sec} is devoted to the proof of Theorem \ref{singularthm}, where we analyze the low and high frequency components separately based on a Littlewood-Paley type decomposition. 
One of the key estimates for the high frequency part in the proof of Theorem \ref{singularthm} is Proposition \ref{mainpropo} and its proof is given in Section \ref{mainpropoprsec}.

\hfill

{\bf Notation}
Let $L^1_{\mathrm{loc}}(\bbrn)$ be the space of all locally integrable functions on $\bbrn$ and $L^{\infty}_c(\bbrn)$ denote the space of all essentially bounded measurable functions with compact support on $\bbrn$.
 We   use the symbol $A\lesssim B$ to indicate that $A\leq CB$ for some constant $C>0$ independent of the variable quantities $A$ and $B$, and $A\sim B$ if $A\lesssim B$ and $B\lesssim A$ hold simultaneously. For each $\ppp:=(p_1,\dots,p_m)$ and $r>0$, we write $r\ppp:= (rp_1,\dots,rp_m)$ and $\ppp/r:=(p_1/r,\dots,p_m/r)$.
For an $m$-tuple of weights $\www=(w_1,\dots,w_m)$ and $\delta>0$, we define $\www^{\delta}:=\big(w_1^{\delta},\dots,w_m^{\delta}\big)$.

 \hfill

\section{Preliminaries}\label{prelisec}

In this section, we provide several auxiliary results that are essential for the proof of Theorem \ref{singularthm}.

\subsection{Maximal inequalities}\label{maxinesec}

 For  a locally integrable function $f$ on $\bbrn$, we  define the  Hardy-Littlewood maximal function by
$$\mathcal{M}f(x):=\sup_{Q:x\in Q}\frac{1}{|Q|}\int_{Q}\big|f(u) \big|\; du $$
where  the supremum is taken over all cubes with sides parallel to the axes containing $x$.
For $0<r<\infty$, we also define the $L^r$-variant of the Hardy-Littlewood maximal operator by the formula
$$\mathcal{M}_rf(x):=\Big(\mathcal{M}\big(|f|^r\big)(x) \Big)^{1/r}.$$
For $0<r<\infty$, the (homogeneous) sharp maximal function $\mathscr{M}_r^{\sharp}f$  is defined by
 $$\mathscr{M}_r^{\sharp}f(x):=\sup_{Q:x\in Q}\inf_{c_Q\in \mathbb{C}}\bigg(\frac{1}{|Q|}\int_Q \big| f(y)-c_Q\big|^r\; dy\bigg)^{1/r}$$
  where the supremum is taken over all cubes in $\bbrn$ containing the point $x$.
  By H\"older's inequality,
\begin{equation}\label{increasingpro}
\mathscr{M}_r^{\sharp}f(x)\le \mathscr{M}_s^{\sharp}f(x)\quad \text{for all }~0<r<s<\infty.
\end{equation}
It is also clear that
$$\mathscr{M}_r^{\sharp}f(x)\lesssim \mathcal{M}_rf(x) \quad \text{ for all }~0<r<\infty$$
and
\begin{equation}\label{m1sharpfrxest}
\Big(\mathscr{M}_1^{\sharp}\big(|f|^r\big)(x)\Big)^{1/r}\le \mathscr{M}_r^{\sharp}f(x)\quad \text{ for all }~0<r\le 1.
\end{equation}
Moreover, if $1\le  p_0\le p<\infty$ and $\mathcal{M}f\in L^{p_0}(\bbrn)$,  then
\begin{equation*}
\Vert \mathcal{M}f\Vert_{L^p(\bbrn)}\lesssim \big\Vert \mathscr{M}^{\sharp}_1f\big\Vert_{L^p(\bbrn)},
\end{equation*} which was established by Fefferman and Stein \cite{Fe_St1972}.
The above inequality also holds in the weighted setting.
\begin{customlemma}{D}\cite[IV. Theorem 2.20]{Ga_Ru1985}\label{sharpmaxweight}
Let $w\in A_{\infty}$, and suppose that $\mathcal{M}f\in L^{p_0}(\bbrn)$ for some $p_0$ with $0<p_0<\infty$. 
Then for every $p_0\le p<\infty,$
we have
$$\big\Vert \mathcal{M}f\big\Vert_{L^p(w)}\lesssim \big\Vert \mathscr{M}^{\sharp}_1f\big\Vert_{L^p(w)}.$$
\end{customlemma}

Combining Lemma \ref{sharpmaxweight}, together with \eqref{increasingpro} and \eqref{m1sharpfrxest},  we have the following result:
For any $0<r<\infty$ and $w\in A_{\infty}$, if $\mathcal{M}_rf\in L^p(\bbrn)$, then
\begin{equation}\label{eq:flpwlemrshflpw}
\Vert f\Vert_{L^p(w)}\lesssim \big\Vert \mathscr{M}_r^{\sharp}f\big\Vert_{L^p(w)}.
\end{equation}
Indeed, setting $s=\min\{1,r\}\le 1$,
\begin{align*}
\Vert f\Vert_{L^p(w)}&\le \big\Vert \mathcal{M}_sf\big\Vert_{L^p(w)}=\Big\Vert \mathcal{M}\big(|f|^s\big)\Big\Vert_{L^{p/s}(w)}^{1/s}\\
&\lesssim \Big\Vert \mathscr{M}_1^{\sharp}\big(|f|^s\big)\Big\Vert_{L^{p/s}(w)}^{1/s}\le \big\Vert \mathscr{M}_s^{\sharp}f\big\Vert_{L^p(w)}\le \big\Vert \mathscr{M}_r^{\sharp}f\big\Vert_{L^p(w)}.
\end{align*}

Let us define
 a multi-sublinear version of Hardy-Littlewood maximal operator $\mathbf{M}$  by
 $$\mathbf{M}\big(f_1,\dots,f_m \big)(x):=\sup_{Q:x\in Q}\bigg( \frac{1}{|Q|^m}\int_{Q^m}\prod_{j=1}^{m}\big|f_j(u_j) \big|\; d\uuu \bigg)$$
 for locally integrable functions $f_1,\dots,f_m$ on $\bbrn$, where $Q^m:=Q\times \cdots \times Q$, $d\uuu:=du_1\cdots du_m$, and
 the supremum is taken over all cubes $Q$ in $\bbrn$ containing $x$.
 Clearly, this maximal function is controlled by the product of maximal functions $\mathcal{M}f_j(x)$.
We also define the $L^r$ version of $\mathbf{M}$ by
 $$\mathbf{M}_r(f_1,\dots,f_m)(x):=\big( \mathbf{M}\big(|f_1|^r,\dots, |f_m|^r \big)(x) \big)^{{1}/{r}}.$$

The following lemma is a very helpful tool to handle multiple weighted estimates, which will be repeatedly used in the proof of Theorem \ref{singularthm}.
\begin{lemma}\label{multmaximalpt}
Let $x\in\bbrn$ and $L>mn$. Then we have
\begin{align*}
\sup_{\gamma>0}\int_{(\bbrn)^m}\frac{1}{\gamma^{mn}}\frac{1}{(1+\frac{|(x-y_1,\dots,x-y_m)|}{\gamma})^{L}}\prod_{j=1}^{m}|f_j(y_j)|\; d\yyy\lesssim_{L,m,n} \mathbf{M}\big(f_1,\dots,f_m\big)(x)
\end{align*}
for all $f_1,\dots,f_m\in L^{\infty}_c(\bbrn)$.
\end{lemma}
\begin{proof}
For any $t>0$,
let $Q(x;t)$ denote the cube, centered at $x$, whose side-length is $t$,
and
\begin{equation}\label{qqxtdef}
\mathbf{Q}(x;t):=\big\{\yyy\in (\bbrn)^m:  y_1,\dots,y_m\in Q(x;t)\big\}.
\end{equation}
Then 
\begin{align*}
&\int_{(\bbrn)^m}\frac{1}{\gamma^{mn}}\frac{1}{(1+\frac{|(x-y_1,\dots,x-y_m)|}{\gamma})^{L}}\prod_{j=1}^{m}|f_j(y_j)|\; d\yyy\\
&=\int_{\mathbf{Q}(x;\gamma)}\frac{1}{\gamma^{mn}}\frac{1}{(1+\frac{|(x-y_1,\dots,x-y_m)|}{\gamma})^{L}}\prod_{j=1}^{m}|f_j(y_j)|\; d\yyy\\
&\qq+\sum_{l=1}^{\infty}\int_{\mathbf{Q}(x;2^{l}\gamma)\setminus \mathbf{Q}(x;2^{l-1}\gamma)}\frac{1}{\gamma^{mn}}\frac{1}{(1+\frac{|(x-y_1,\dots,x-y_m)|}{\gamma})^{L}}\prod_{j=1}^{m}|f_j(y_j)|\; d\yyy\\
&\lesssim \sum_{l=0}^{\infty}2^{-l(L-mn)}\frac{1}{(2^l\gamma)^{mn}}\prod_{j=1}^{m}\int_{Q(x;2^l\gamma)}\big| f_j(y_j)\big|\; dy_j\\
&\lesssim_L \mathbf{M}\big(f_1,\dots,f_m\big)(x).
\end{align*}
The desired result is immediate by taking the supremum over $\gamma>0$.
\end{proof}

\hfill

\subsection{Multiple weight classes}

The following lemma gives a structural decomposition of the multiple weight class $\mathrm{A}_{\ppp}$ into componentwise Muckenhoupt $A_p$ classes.
\begin{customlemma}{E}\cite[Theorem 3.6]{Le_Om_Pe_To_Tr2009}\label{multiindiweight}
Let $\www=(w_1,\dots,w_m)$ and $1\le p_1,\dots,p_m<\infty$ with $1/p=1/p_1+\cdots+1/p_m$. Then
$$\www\in \mathrm{A}_{\ppp}$$
if and only if 
$$\begin{cases}
w_j^{1-p_j'}\in A_{mp_j'}, & j=1,\dots,m,\\
v_{\www}\in A_{mp},&
\end{cases}$$
where the condition $w_j^{1-p_j'}\in A_{mp_j'}$ in the case $p_j=1$ is understood as $w_j^{1/m}\in A_1$.
\end{customlemma}
This result enables us to extend  properties of $A_p$ weight to the multilinear setting, considering each component separately. Indeed, this lemma serves as a backbone for the proof of Lemmas \ref{weightlowpore}, \ref{mwextra}, and \ref{mweightincl}.

It follows from H\"older's inequality that for $1\le p<\infty$
\begin{equation}\label{linearincl}
w\in A_p\qq \Rightarrow \qq  w^{\delta}\in A_p, \q 0\le \delta\le 1.
\end{equation}
This property also extends to multiple weight spaces $\mathrm{A}_{\ppp}$, owing to Lemma \ref{multiindiweight}.
\begin{lemma}\label{weightlowpore}
Let $1\le p_1,\dots,p_m<\infty$ with $1/p=1/p_1+\cdots+1/p_m$. Suppose that $\www\in \mathrm{A}_{\ppp}$.
Then for any $0\le \delta\le 1$,
$$\www^{\delta}\in \mathrm{A}_{\ppp}.$$
\end{lemma}
\begin{proof}
It is trivial when $\delta=0$ or $\delta=1$.
According to Lemma \ref{multiindiweight}, we have
$$\begin{cases}
w_j^{1-p_j'}\in A_{mp_j'}, & j=1,\dots,m,\\
v_{\www}\in A_{mp},&
\end{cases}$$
and then \eqref{linearincl} yields that
$$\begin{cases}
(w_j^{\delta})^{1-p_j'}\in A_{mp_j'}, & j=1,\dots,m,\\
v_{\www}^{\delta}\in A_{mp},&
\end{cases}$$
which is equivalent to $\www^{\delta}\in \mathrm{A}_{\ppp}$, using Lemma \ref{multiindiweight} again.
\end{proof}

We also recall from \cite[IV. Theorem 2.7]{Ga_Ru1985} that when $1\le p<\infty$, 
 \begin{equation}\label{w1eapcon}
\text{ $w\in A_p~$ implies $~w^{1+\epsilon}\in A_p\q$ for some $\epsilon>0$}.
 \end{equation}
An analogous property also holds in the setting of multiple weight classes.
\begin{lemma}\label{mwextra}
Let $1\le p_1,\dots,p_m<\infty$.  Assume that $\www\in \mathrm{A}_{\ppp}$. Then there exists $\epsilon>0$ such that
$$\www^{1+\epsilon}\in \mathrm{A}_{\ppp}.$$
\end{lemma}
\begin{proof}
Assume that $\www\in \mathrm{A}_{\ppp}$ and define $p$ by $1/p=1/p_1+\cdots+1/p_m$.
By Lemma \ref{multiindiweight}, we have
$$v_{\www}\in A_{mp}, \qq w_j^{1-p_j'}\in A_{mp_j'},~ j=1,\dots,m.$$
We apply \eqref{w1eapcon} to choose  $\epsilon_0,\epsilon_1,\dots,\epsilon_m>0$ such that
$$v_{\www}^{1+\epsilon_0}\in A_{mp}, \qq w_j^{(1+\epsilon_j)(1-p_j')}\in A_{mp_j'}, ~ j=1,\dots,m.$$
Let 
$$\epsilon:=\min \{\epsilon_0,\epsilon_1,\dots,\epsilon_m\}>0$$
and then \eqref{linearincl} yields
$$v_{\www}^{1+\epsilon}\in A_{mp}, \qq w_j^{(1+\epsilon)(1-p_j')}\in A_{mp_j'}, ~ j=1,\dots,m,$$
which is equivalent to 
$$\www^{1+\epsilon}\in \mathrm{A}_{\ppp},$$
as desired.
\end{proof}

As discussed in \cite[Remark 7.3]{Le_Om_Pe_To_Tr2009}, the classes $\mathrm{A}_{\ppp}$ do not satisfy an increasing inclusion property, unlike \eqref{apweightinclusion}.
That is, even though $p_j\le q_j$ for all $j=1,\dots,m$, the inclusion $\mathrm{A}_{\ppp}\subset A_{\qqq}((\bbrn)^m)$ is not guaranteed in general. However, an inclusion relation holds if the ratios $p_j/q_j$ are fixed for all $j=1,\dots,m$.
\begin{lemma}\label{mweightincl}
Let $1\le p_1,\dots,p_m<\infty$ and $r>1$.
Then we have
$$\mathrm{A}_{\ppp} \subset \mathrm{A}_{r\ppp}.$$
\end{lemma}
\begin{proof}
Assume that $\www\in \mathrm{A}_{\ppp}$ and $r>1$, and define $p$ by $1/p=1/p_1+\cdots+1/p_m$.
By Lemma \ref{multiindiweight}, we have
$$v_{\www}\in A_{mp}, \qq w_j^{1-p_j'}\in A_{mp_j'},~ j=1,\dots,m.$$
First of all, \eqref{apweightinclusion} implies
$$v_{\www}\in A_{mrp}.$$
Moreover, if $p_j>1$,  we have
\begin{align*}
&\Big(\frac{1}{|Q|}\int_Q w_j(x)^{1-(rp_j)'}\; dx \Big) \Big(\frac{1}{|Q|}\int_Q w_j(x)^{( 1-(rp_j)' )(-\frac{1}{m(rp_j)'-1})} \, dx\Big)^{m(rp_j)'-1}\\
&= \Big(\frac{1}{|Q|}\int_Q w_j(x)^{-\frac{1}{rp_j-1}}\; dx \Big)\Big( \frac{1}{|Q|}\int_Q w_j(x)^{\frac{1}{(m-1)rp_j+1}}\; dx\Big)^{\frac{(m-1)rp_j+1}{rp_j-1}}\\
&\le \bigg[\Big(\frac{1}{|Q|}\int_Q w_j(x)^{-\frac{1}{p_j-1}}\; dx \Big)\Big(\frac{1}{|Q|}\int_Q w_j(x)^{\frac{1}{(m-1)p_j+1}} \; dx\Big)^{\frac{(m-1)p_j+1}{p_j-1}}\bigg]^{\frac{p_j-1}{rp_j-1}}\\
&=\bigg[ \Big(\frac{1}{|Q|}\int_Q w_j(x)^{1-p_j'}\; dx \Big) \Big(\frac{1}{|Q|}\int_Q w_j^{( 1-p_j' )(-\frac{1}{m p_j'-1})}\, dx \Big)^{mp_j'-1} \bigg]^{\frac{p_j-1}{rp_j-1}}
\end{align*}
where the inequality follows simply from H\"older's inequality with $\frac{rp_j-1}{p_j-1}>1$ and $\frac{(m-1)rp_j+1}{(m-1)p_j+1}>1$.
This yields that
$$ [w_j^{1-(rp_j)'}]_{A_{m(rp_j)'}} \le [w_j^{1-p_j'}]_{A_{mp_j'}}^{\frac{p_j-1}{rp_j-1}}<\infty, \q j=1,\dots,m$$
and thus
$$w_j^{1-(rp_j)'}\in A_{m(rp_j)'}, \q j=1,\dots,m.$$
If $p_j=1$, then Lemma \ref{multiindiweight} gives $w_j^{1/m}\in A_1$.
By the standard endpoint modification, this implies
$w_j^{1-r'}\in A_{mr'}$, which is precisely the required component condition
corresponding to the exponent $rp_j=r$.

Now Lemma \ref{multiindiweight} concludes
$\www\in \mathrm{A}_{r\ppp},$
as desired.
\end{proof}

\medskip

We shall also use the openness property of Muckenhoupt weights. It is well
known that if $1<p<\infty$ and $w\in A_p$, then there
exists $1<q<p$ such that
\begin{equation}\label{apopen}
w\in A_q.
\end{equation}
See, for example, \cite[IV. Theorem 2.6]{Ga_Ru1985}. The following
multiple-weight version will be used in the proof of Theorem
\ref{singularthm}.

\begin{lemma}\label{mweightopen}
Let $1<p_1,\dots,p_m<\infty$.
Assume that
$\www\in \mathrm{A}_{\ppp}.$
Then there exists $s>1$ such that
$$\www\in \mathrm{A}_{\ppp/s},$$
where
$\ppp/s:=({p_1}/{s},\dots,{p_m}/{s}).$
\end{lemma}

\begin{proof}
Define $p$ by $1/p=1/p_1+\cdots+1/p_m$.
By Lemma \ref{multiindiweight}, we have
$$v_{\www}\in A_{mp} \qquad \text{ and }\qquad w_j^{1-p_j'}\in A_{mp_j'}, \quad j=1,\dots,m.$$
By the scalar openness property \eqref{apopen}, there exist numbers
$$1<a<mp \qquad \text{ and }\qquad 1<a_j<mp_j'=\frac{mp_j}{p_j-1}, \quad  j=1,\dots,m,$$
such that
$$v_{\www}\in A_a \quad \text{ and } \quad w_j^{1-p_j'}\in A_{a_j},~ j=1,\dots,m.$$
By the self-improvement property \eqref{w1eapcon}, for each $j=1,\dots,m$, there exists $\eta_j>0$ such that
$$\big(w_j^{1-p_j'}\big)^{1+\eta_j}\in A_{a_j}.$$
Now we take $s>1$ sufficiently close to $1$ so that
$$\frac{mp}{s}>a, \quad \frac{mp_j}{p_j-s}>a_j,\quad   \frac{p_j}{s}>1,\quad \text{and}\quad 1<\frac{s(p_j-1)}{p_j-s}<1+\eta_j, \qquad j=1,\dots,m.$$
Then, by the monotonicity of scalar Muckenhoupt classes,
$$v_{\www}\in A_{mp/s}.$$
Set
$$\widetilde p_j:=\frac{p_j}{s},\qquad j=1,\dots,m, \qquad\text{and}\qquad \frac1{\widetilde p}=\frac1{\widetilde p_1}+\cdots+\frac1{\widetilde p_m}.$$
Then $\widetilde p=p/s$, and
\begin{equation}\label{wtvwwwinamwtp}
\widetilde v_{\www}:=\prod_{j=1}^m w_j^{\widetilde p/\widetilde p_j}=\prod_{j=1}^m w_j^{p/p_j}=v_{\www}\in A_{m\widetilde p}.
\end{equation}
Moreover, since $\widetilde{p_j}'=\frac{p_j}{p_j-s}$,
we have
$$w_j^{1-\widetilde{p_j}'}=w_j^{-\frac{s}{p_j-s}}=\big(w_j^{1-p_j'}\big)^{\frac{s(p_j-1)}{p_j-s}}.$$
Since $\frac{s(p_j-1)}{p_j-s}<1+\eta_j$,
the property \eqref{linearincl} applied to $\big(w_j^{1-p_j'}\big)^{1+\eta_j}\in A_{a_j}$ gives
$$w_j^{1-\widetilde{p_j}'}=\big(w_j^{1-p_j'}\big)^{\frac{s(p_j-1)}{p_j-s}}\in A_{a_j}\subset A_{m\wt{p_j}'}, \quad j=1,\dots,m.$$
Together with \eqref{wtvwwwinamwtp},
Lemma \ref{multiindiweight} yields
$$\www\in \mathrm{A}_{\wt{{\ppp}}}=\mathrm{A}_{\ppp/s}$$
where $\wt{\ppp}:=(\wt{p_1},\dots,\wt{p_m})$.
This completes the proof.
\end{proof}

\medskip

Recall from \cite{Mu1972} that
 for $1<p<\infty$ 
  \begin{equation}\label{chaweight}
  w\in A_p \q\text{if and only if}\q \Vert \mathcal{M}f\Vert_{L^p(w)}\lesssim \Vert f\Vert_{L^p(w)}.
  \end{equation}
As a multilinear analogue of  \eqref{chaweight}, 
the class $\mathrm{A}_{\ppp}$ can be also characterized by a maximal inequality. 
   \begin{customlemma}{F}\cite{Le_Om_Pe_To_Tr2009}\label{multiweightcha}
Let $1<p_1,\dots,p_m<\infty$ with $1/p=1/p_1+\dots+1/p_m$.
Then the inequality
 $$\big\Vert \mathbf{M}(f_1,\dots,f_m)\big\Vert_{L^p(v_{\www})}\lesssim \prod_{j=1}^{m}\Vert f_j\Vert_{L^{p_j}(w_j)}$$
 holds for all locally integrable functions $f_1,\dots,f_m$ if and only if
 $$\www=(w_1,\dots,w_m)\in \mathrm{A}_{\ppp}.$$
 \end{customlemma}
As a consequence of Lemmas \ref{multiweightcha} and \ref{mweightincl}, we can obtain that when $1<r<p_1,\dots,p_m<\infty$, 
\begin{equation}\label{mrweightest}
\big\Vert \mathbf{M}_r(f_1,\dots,f_m)\big\Vert_{L^p(v_{\www})}\lesssim \prod_{j=1}^{m}\Vert f_j\Vert_{L^{p_j}(w_j)},
\end{equation}
provided that $\www\in \mathrm{A}_{\ppp/r}$.

\hfill

\subsection{Multilinear (convolution-type) Calder\'on-Zygmund operators}

Let $K$ be a locally integrable function defined on $(\bbrn)^{m}\setminus \{\000\}$ satisfying the size estimate
\begin{equation}\label{multiczsize}
\big| K(y_1,\dots,y_m)\big|\le \frac{A}{|(y_1,\dots,y_m)|^{mn}}
\end{equation}
and the smoothness estimate
\begin{equation}\label{multiczsm}
\big|K(y_1,\dots,y_j,\dots,y_m)-K(y_1,\dots,y_j',\dots,y_m) \big|\le \frac{A|y_j-y_j'|^{\epsilon_0}}{|(y_1,\dots,y_m)|^{mn+\epsilon_0}}
\end{equation}
for some $\epsilon_0>0$ and all $1\le j\le m$, whenever $2|y_j-y_j'|\le \max_{1\le k\le m}|y_k|$.
Then the associated $m$-linear singular integral operator $T$ is defined by
$$T\big(f_1,\dots,f_m\big)(x):=\mathrm{p.v.} \int_{(\bbrn)^m}K(\yyy)\prod_{j=1}^{m}f_j(x-y_j)\; d\yyy$$
for $f_1,\dots,f_m\in \mathscr{C}_c^{\infty}(\bbrn)$. 
We say that $T$ is an $m$-linear Calder\'on-Zygmund operator when it satisfies
$$\big\Vert T(f_1,\dots,f_m)\big\Vert_{L^q(\bbrn)}\le B \prod_{j=1}^{m}\Vert f_j\Vert_{L^{q_j}(\bbrn)}$$
for some $1<q_1,\dots,q_m<\infty$ with $1/q=1/q_1+\cdots+1/q_m$, and for some $B>0$.
Then it is known in \cite{Gr_To2002} that 
\begin{equation}\label{strongtf1mest}
\big\Vert T(f_1,\dots,f_m)\big\Vert_{L^p(\bbrn)}\lesssim_{\ppp} (A+B)\prod_{j=1}^{m}\Vert f_j\Vert_{L^{p_j}(\bbrn)}
\end{equation}
when $1<p_1,\dots,p_m<\infty$ with $1/p=1/p_1+\cdots+1/p_m$.

It is also known in \cite{GT-Indiana} that if $T$ is an $m$-linear Calder\'on-Zygmund operator, then 
\begin{equation}\label{ptdefalev}
\text{$T(f_1,\dots,f_m)$ is pointwise well-defined almost everywhere}
\end{equation}
 when $f_j\in L^{p_j}(\bbrn)$ with $1< p_j< \infty$, and 
 \begin{equation}\label{boundedforlpj}
 \text{ \eqref{strongtf1mest} also works   for  $f_j\in L^{p_j}(\bbrn)$ with $1< p_j< \infty$.}
 \end{equation}

\hfill

\subsection{Complex Interpolation for multilinear operators}

Let $\mathbf{S}:=\big\{z\in\mathbb{C}: 0< \mathrm{Re}z < 1 \big\}$ be the open unit strip on the complex plane $\mathbb{C}$ and $\overline{\mathbf{S}}$ be its closure.
Suppose that for every $z\in \overline{\mathbf{S}}$, $T_z$ is an $m$-linear operator defined on $\mathscr{C}_c^{\infty}(\bbrn)\times\cdots\times \mathscr{C}_c^{\infty}(\bbrn)$ taking values in $L_{\mathrm{loc}}^1(\bbrn)$. We call $\{T_z\}_{z}$ an analytic family if for all $f_1,\dots,f_m\in \mathscr{C}_c^{\infty}(\bbrn)$ and $w$ bounded function with compact support on $\bbrn$ the mapping
$$z \mapsto \int_{\bbrn}T_z\big(f_1,\dots,f_m\big)(x)w(x)\; dx$$
is analytic in the open strip $\mathbf{S}$ and continuous on $\overline{\mathbf{S}}$.
The analytic family $\{T_z\}_z$ is called of admissible growth if there is a constant $\gamma$ with $0\le \gamma<\pi$ and $1\le t\le \infty$ such that for any $f_1,\dots,f_m\in\mathscr{C}_c^{\infty}(\bbrn)$ and every compact set $\mathcal{B}$ in $\bbrn$, there exists a constant $C_{f_1,\dots,f_m,\mathcal{B}}$ such that
\begin{equation*}
\log{\bigg(\int_B\big| T_z\big(f_1,\dots,f_m\big)(x)\big|^t\; dx \bigg)^{1/t}}\le C_{f_1,\dots,f_m,\mathcal{B}}e^{\gamma |\mathrm{Im}(z)|}, \q \text{ for all }~ z\in \overline{\mathbf{S}}.
\end{equation*}

Then we have the following interpolation theory for analytic multilinear operators, which is a multilinear version of Stein's complex interpolation for analytic families in \cite{St1956}.
\begin{customlemma}{G}\cite{Gr_Ou2022}\label{weightinterpol}
For $z\in \overline{\mathbf{S}}$, let $T_{z}$ be an $m$-linear operator on $\mathscr{C}_c^{\infty}(\bbrn)\times \cdots\times \mathscr{C}_c^{\infty}(\bbrn)$ with values in $L_{\mathrm{loc}}^{1}(\bbrn)$ that form an analytic family of admissible growth. For $j=1,\dots,m$ let $0<p_j^0,p^1_j< \infty$, $0<p^0,p^1< \infty$, and suppose that $0<p,p_1,\dots,p_m< \infty$ satisfy
$$\frac{1}{p_j}=\frac{1-\theta}{p_j^0}+\frac{\theta}{p_j^1}\q \text{ and }\q \frac{1}{p}=\frac{1-\theta}{p^0}+\frac{\theta}{p^1}$$
for some $0<\theta<1$.
Suppose that for all $f_1,\dots,f_m\in\mathscr{C}_c^{\infty}(\bbrn)$ we have
\begin{equation}\label{iniconininter}
\begin{aligned}
\big\Vert T_{iy}(f_1,\dots,f_m)\big\Vert_{L^{p^0}(\bbrn)}&\le M_0 \prod_{j=1}^{m}\Vert f_j\Vert_{L^{p_j^0}(\bbrn)}\\
\big\Vert T_{1+iy}(f_1,\dots,f_m)\big\Vert_{L^{p^1}(\bbrn)}&\le M_1 \prod_{j=1}^{m}\Vert f_j\Vert_{L^{p_j^1}(\bbrn)}.
\end{aligned}
\end{equation}
Then  we have
$$\big\Vert T_{\theta}(f_1,\dots,f_m)\big\Vert_{L^p(\bbrn)}\le M_0^{1-\theta}M_1^{\theta}\prod_{j=1}^{m}\Vert f_j \Vert_{L^{p_j}(\bbrn)}$$
for $f_1,\dots,f_m\in \mathscr{C}_c^{\infty}(\bbrn)$.
\end{customlemma}
We remark that the original version in \cite[Theorem 3.2]{Gr_Ou2022} provides more general cases in the setting that functions $f_1,\dots,f_m$ are defined on metric measure spaces in which balls have finite measure, and the bounds $M_0$ and $M_1$ in \eqref{iniconininter} are continuous functions of $y$. Moreover, the original one deals with continuous functions $f_j$ with compact support, which originated from Urysohn's lemma, stated in \cite[Lemma 2.1]{Gr_Ou2022}. However,  all of the arguments therein can be also valid with $f_j\in \mathscr{C}_c^{\infty}(\bbrn)$, simply applying a smooth version of Urysohn's lemma (see \cite[page 38]{Li_Lo2001}), which allows constructions of $f_{z}^{\epsilon}\in \mathscr{C}_c^{\infty}(\bbrn)$ in  \cite[(2.1)]{Gr_Ou2022}.

\hfill

\section{Proof of Theorem \ref{singularthm}}\label{prthm1sec}

When $q=\infty$, and the weight assumption becomes $\www\in \mathrm{A}_{\ppp}.$
By Lemma \ref{mweightopen}, there exists $s>1$ such that $\www\in \mathrm{A}_{\ppp/s}$.
Choose $q_0<\infty$ sufficiently large so that
$$ 1<q_0'<\min\{s,p_1,\dots,p_m\}.$$
Then Lemma \ref{mweightincl} gives
$ \mathrm{A}_{\ppp/s} \subset \mathrm{A}_{\ppp/q_0'}$ and thus
$\www\in \mathrm{A}_{\ppp/q_0'}.$
Moreover, $\Omega\in L^\infty(\mathbb S^{mn-1})\subset L^{q_0}(\mathbb S^{mn-1})$.
Therefore, the desired estimate in the case $q=\infty$ follows from the
finite $q$ case applied with $q=q_0$. Hence, throughout the proof, it is
enough to assume $1<q<\infty$.

Now let us assume $1<q<\infty$. 
Let $\Psi$ be a Schwartz function on $(\bbrn)^m$ whose Fourier transform is supported in the annulus $\{\xxxi\in (\bbrn)^m: \frac{1}{2}\le |\xxxi|\le 2\}$ and satisfies
$\sum_{k\in\bbz}\wh{\Psi_k}(\xxxi)=1$ for $\xxxi\not= \000$ where $\wh{\Psi_k}(\xxxi):=\wh{\Psi}(2^{-k} \xxxi )$.
For each $\gamma \in\bbz$ and $\mu\in\bbz$, let
 \begin{equation*}
  K^{\gamma}(\yyy):=\wh{\Psi}(2^{\gamma}\yyy)K(\yyy), \q   K_{\mu}^{\gamma}(\yyy):=\Psi_{{\mu}+\gamma}\ast K^{\gamma}(\yyy), \qquad \yyy\in (\bbrn)^m
  \end{equation*}
and 
   \begin{equation}\label{kgmdef}
 K_{\mu}:=\sum_{\gamma\in\bbz}K_{\mu}^{\gamma}.
  \end{equation}
Then the corresponding operators $T_{K_{\mu}^{\gamma}}$ and $T_{K_{\mu}}$  are defined as
\begin{align*}
T_{K_{\mu}^{\gamma}}\big(f_1,\dots,f_m \big)(x)&:=\int_{(\bbrn)^m}K_{\mu}^{\gamma}(\yyy)\prod_{j=1}^{m}f_j(x-y_j)\; d\yyy,\\
T_{K_{\mu}}\big(f_1,\dots,f_m\big)(x)&:=\int_{(\bbrn)^m}K_{\mu}(\yyy)\prod_{j=1}^{m}f_j(x-y_j) \; d\yyy
\end{align*}
  so that
  $$\LL_{\Omega}\big(f_1,\dots,f_m\big)=\sum_{\mu\in\bbz}\sum_{\gamma\in\bbz}T_{K_{\mu}^{\gamma}}\big(f_1,\dots,f_m\big)=\sum_{\mu\in\bbz}T_{K_{\mu}}\big(f_1,\dots,f_m\big).$$
We remark that
 \begin{equation*}
 K^\gamma(\yyy)=2^{\ga mn} K^0(2^\gamma \yyy),
 \end{equation*} 
which deduces
 \begin{equation*}
K_{\mu}^{\gamma}(\yyy)=2^{\gamma mn}\big(\Psi_{{\mu}}\ast K^{0}\big)(2^\gamma \yyy)=2^{\gamma mn}K_{\mu}^{0}(2^{\gamma}\yyy)
\end{equation*}
or equivalently,
$$
\wh{K^\gamma_\mu}(\xxxi)= \wh{\Psi_{\mu+\gamma}}(\xxxi)\wh{K^0}(2^{-\ga}\xxxi)=\wh{K^0_\mu}(2^{-\ga} \xxxi).
$$

Duoandikoetxea and Rubio de Francia \cite{Du_Ru1986} proved that if $0<\delta<\frac{1}{q'}$, then 
 \begin{align*}
 \big| \wh{K^0}(\xxxi)\big|&\lesssim \Vert \Omega\Vert_{L^q(\mathbb{S}^{mn-1})}\min\big\{|\xxxi|,|\xxxi|^{-\delta} \big\}\\
  \big| \partial^{\alpha}\wh{K^0}(\xxxi)\big|&\lesssim \Vert \Omega\Vert_{L^q(\mathbb{S}^{mn-1})}\min\big\{1,|\xxxi|^{-\delta} \big\}, \qq \alpha\not= \000
 \end{align*}
 and accordingly,
  \begin{align*}
 \Big| \wh{K_{\mu}}(\xxxi)\Big|&\lesssim \Vert \Omega\Vert_{L^q(\mathbb{S}^{mn-1})}\min\big\{2^{\mu},2^{-\delta \mu} \big\}\\
  \Big| \partial^{\alpha}\wh{K_{\mu}}(\xxxi)\Big|&\lesssim \Vert \Omega\Vert_{L^q(\mathbb{S}^{mn-1})} \mathcal{G}(\mu,\alpha)\frac{1}{|\xxxi|^{|\alpha|}}, \qq 1\le |\alpha| 
 \end{align*}
 where
 $$\mathcal{G}(\mu,\alpha):=\begin{cases}
 2^{\mu(|\alpha|-\delta)}& \mu>0\\
 2^{\mu}& \mu\le 0
 \end{cases}.$$
 See the proof of \cite[Proposition 3]{Gr_He_Ho2018} for more details.
Then a standard computation yields that
\begin{equation}\label{kmuyyyptest}
\big| K_{\mu}(\yyy)\big|\lesssim \Vert \Omega\Vert_{L^q(\mathbb{S}^{mn-1})}\frac{1}{|\yyy|^{mn}}\begin{cases}
2^{\mu(mn+1-\delta)}& \mu>0\\
2^{\mu} &\mu\le 0
\end{cases}
\end{equation}
and
\begin{equation}\label{nakmuyyyptest}
\big| \nabla K_{\mu}(\yyy)\big|\lesssim \Vert \Omega\Vert_{L^q(\mathbb{S}^{mn-1})}\frac{1}{|\yyy|^{mn+1}}\begin{cases}
2^{\mu (mn+2-\delta)} & \mu>0\\
2^{\mu} & \mu\le 0
\end{cases}.
\end{equation}

\hfill

\subsection{The case when $\mu\le 0$}

We first set
\begin{equation*}
\mathcal{K}:=\sum_{\mu\in\bbz:\mu\le 0}K_{\mu}.
\end{equation*}
Then
\begin{align*}
\sum_{\mu\in\bbz:\mu\le 0}T_{K_{\mu}}\big(f_1,\dots,f_m\big)(x)&=\int_{(\bbrn)^m}\mathcal{K}(\yyy) \prod_{j=1}^{m}f_j(x-y_j) \; d\yyy\\
&=:T_{\mathcal{K}}\big(f_1,\dots,f_m\big)(x).
\end{align*}
The (unweighted) $L^{p_1}\times \cdots \times L^{p_m}\to L^p$ boundedness for $T_{\mathcal{K}}$, $1<p_1,\dots,p_m<\infty$, was already verified in \cite{Park_submitted2} with the constant $C\Vert \Omega\Vert_{L^q(\mathbb{S}^{mn-1})}$.
Moreover, the estimates \eqref{kmuyyyptest} and \eqref{nakmuyyyptest} deduce that the kernel $\mathcal{K}$ satisfies the size and smoothness conditions \eqref{multiczsize} and \eqref{multiczsm} for multilinear Calder\'on-Zygmund kernel with constant $C\Vert \Omega\Vert_{L^q(\mathbb{S}^{mn-1})}$. 
Now it follows from \cite[Corollary 3.9]{Le_Om_Pe_To_Tr2009}  that for  $1<p_1,\dots,p_m<\infty$ with $1/p=1/p_1+\cdots+1/p_m,$
\begin{equation}\label{tkf1fmlpvwest}
\big\Vert T_{\mathcal{K}}(f_1,\dots,f_m)\big\Vert_{L^p(v_{\www})}\lesssim \Vert \Omega\Vert_{L^q(\mathbb{S}^{mn-1})}\prod_{j=1}^{m}\Vert f_j\Vert_{L^{p_j}(w_j)},
\end{equation}
provided that $\www\in \mathrm{A}_{\ppp}$, where $f_j\in L^{\infty}_c(\bbrn)\cap L^{p_j}(w_j)$, $1\le j\le m$. 
Here, we recall that $\mathrm{A}_{\ppp/q'}\subset \mathrm{A}_{\ppp}$, in view of Lemma \ref{mweightincl}.

\hfill

\subsection{The case when $\mu>0$}

We first recall from \cite[Claim 5]{Do_Sl2024}  that for any $1<p_1,\dots,p_m<\infty$ with $1/p=1/p_1+\dots+1/p_m$ and $1<q<\infty$ satisfying \eqref{qcondequi}, there exists $\delta_0>0$ such that
\begin{equation}\label{nonweightkeyest}
\big\Vert T_{K_{\mu}}(f_1,\dots,f_m)\big\Vert_{L^p(\bbrn)}\lesssim 2^{-\delta_0 \mu}\Vert \Omega\Vert_{L^q(\mathbb{S}^{mn-1})}\prod_{j=1}^{m}\Vert f_j\Vert_{L^{p_j}(\bbrn)}
\end{equation}
for any $f_1,\dots,f_m\in\mathscr{C}_c^{\infty}(\bbrn)$.

Moreover, the multilinear singular integral operator $T_{K_{\mu}}$ is an $m$-linear Calder\'on-Zygmund operator with constant $2^{\mu(mn+2-\delta)}\Vert \Omega\Vert_{L^q(\mathbb{S}^{mn-1})}$ for $\mu>0$, in view of  \eqref{kmuyyyptest}, \eqref{nakmuyyyptest}, and \eqref{nonweightkeyest}. Therefore, according to \eqref{ptdefalev},
\begin{equation}\label{tkmug1gmptest}
\text{$T_{K_{\mu}}(g_1,\dots,g_m)$ is well-defined almost everywhere when $g_1,\dots,g_m\in L^{\infty}_c(\bbrn)$.}
\end{equation}
Since $L^{\infty}_c(\bbrn)\subset L^{p_j}(\bbrn)$ and $\mathscr{C}_c^{\infty}(\bbrn)$ is dense in $L^{p_j}(\bbrn)$, applying a standard argument to \eqref{nonweightkeyest}, we obtain
\begin{equation}\label{bcrnlpest}
\big\Vert T_{K_{\mu}}(g_1,\dots,g_m)\big\Vert_{L^p(\bbrn)}\lesssim 2^{-\delta_0 \mu}\Vert \Omega\Vert_{L^q(\mathbb{S}^{mn-1})}\prod_{j=1}^{m}\Vert g_j\Vert_{L^{p_j}(\bbrn)}
\end{equation}
for all $g_1,\dots,g_m\in L^{\infty}_c(\bbrn)$.
Moreover, similar to \eqref{tkf1fmlpvwest}, \cite[Corollary 3.9]{Le_Om_Pe_To_Tr2009} yields that if $\www\in \mathrm{A}_{\ppp}$, each $T_{K_{\mu}}$, $\mu>0$, satisfies the weighted estimates
\begin{equation}\label{tkmug1gm2mumn2d}
\big\Vert T_{K_{\mu}}(g_1,\dots,g_m)\big\Vert_{L^p(v_{\www})}\lesssim 2^{\mu(mn+2-\delta)}\Vert \Omega\Vert_{L^q(\mathbb{S}^{mn-1})}\prod_{j=1}^{m}\Vert g_j\Vert_{L^{p_j}(w_j)}
\end{equation}
for all $g_j \in L_c^{\infty}(\bbrn)$, $1\le j\le m$.

\begin{proposition}\label{mainpropo}
Let $1<q< \infty$ and $\mu>0$. Suppose that $r=q'/m$.
Then we have
\begin{equation*}
\mathscr{M}_r^{\sharp}\big(T_{K_{\mu}}(g_1,\dots,g_m)  \big)(x)\lesssim_{q,r} \mu\Vert \Omega\Vert_{L^q(\mathbb{S}^{mn-1})} \mathbf{M}_{q'}\big(g_1,\dots,g_m\big)(x)
\end{equation*}
for $g_1,\dots,g_m\in L^{\infty}_c(\bbrn)$.
\end{proposition}
The proposition will be proved in the next section.
As a consequence of Proposition \ref{mainpropo}, the following weighted norm inequality holds.
\begin{corollary}\label{qconjulessp}
Let $1<\wt{p_1},\dots,\wt{p_m}<\infty$ with $1/\wt{p}=1/\wt{p_1}+\cdots+1/\wt{p_m}$. Suppose that $1<q< \infty$ satisfies $q'< \wt{p_1},\dots,\wt{p_m}$ and $\www\in \mathrm{A}_{\vec{\wt{\boldsymbol{p}}}/q'}$,
 where $\vec{\wt{\boldsymbol{p}}}:=(\wt{p_1},\dots,\wt{p_m})$.
Then for any $\mu>0$, we have
\begin{equation*}
\big\Vert T_{K_{\mu}}(g_1,\dots,g_m)\big\Vert_{L^{\wt{p}}(v_{\www})}\lesssim \mu\Vert \Omega\Vert_{L^q(\mathbb{S}^{mn-1})}\prod_{j=1}^{m}\Vert g_j\Vert_{L^{\wt{p_j}}(w_j)}
\end{equation*}
where $g_j \in L^{\infty}_c(\bbrn)$, $1\le j\le m$.
\end{corollary}
\begin{proof}
Let $r=q'/m$. Then we have
 $1/\wt{p}\le m/\min\{\wt{p_1},\dots,\wt{p_m}\}<1/r$ so that $r<\wt{p}$. 
For any  $g_1,\dots,g_m\in L^{\infty}_c(\bbrn)$, \eqref{bcrnlpest} proves $T_{K_{\mu}}\big(g_1,\dots,g_m\big)\in L^{\wt{p}}(\bbrn)$. 
This implies
$$\mathcal{M}_r\Big( T_{K_{\mu}}(g_1,\dots,g_m) \Big)\in L^{\wt{p}}(\bbrn),$$
in view of the $L^{\wt{p}}$ boundedness of $\mathcal{M}_r$ for $r<\wt{p}$.
Now \eqref{eq:flpwlemrshflpw}, Proposition \ref{mainpropo}, and \eqref{mrweightest} prove
\begin{align*}
\big\Vert T_{K_{\mu}}(g_1,\dots,g_m)\big\Vert_{L^{\wt{p}}(v_{\www})}&\lesssim_r \big\Vert \mathscr{M}_r^{\sharp}\big(T_{K_{\mu}}(g_1,\dots,g_m)  \big)\big\Vert_{L^{\wt{p}}(v_{\www})}\\
&\lesssim \mu\Vert \Omega\Vert_{L^q(\mathbb{S}^{mn-1})}\big\Vert    \mathbf{M}_{q'}(g_1,\dots,g_m)  \big\Vert_{L^{\wt{p}}(v_{\www})}\\
&\lesssim \mu\Vert \Omega\Vert_{L^q(\mathbb{S}^{mn-1})}\prod_{j=1}^{m}\Vert g_j\Vert_{L^{\wt{p_j}}(w_j)}
\end{align*}
as $v_{\www}\in A_{\infty}$.
\end{proof}

Now let us complete the proof of Theorem \ref{singularthm}.
For $q'\le p_1,\dots,p_m$, we assume $\www\in \mathrm{A}_{\ppp/q'}$. 
According to Lemma \ref{mwextra}, there exists $\epsilon>0$ such that
$$\www^{1+\epsilon}\in \mathrm{A}_{\ppp/q'}.$$
We choose $r_0<1$ sufficiently close to $1$ such that
$$\max\Big\{ \frac{\epsilon}{1+\epsilon},\frac{1}{p_1},\dots,\frac{1}{p_m}\Big\}<r_0<1$$
and
$$\sum_{j\in J}\frac{1}{p_j r_0}<\frac{|J|}{q'}+\frac{1}{q} \quad \text{ for every $J\subset \{1,\dots,m\}$}.$$
Setting
$$r_1=\frac{r_0}{r_0-\epsilon(1-r_0)}>1,$$
we have
$$\frac{\frac{\epsilon}{1+\epsilon}}{r_0}+\frac{\frac{1}{1+\epsilon}}{r_1}=1,$$
which implies
$$\frac{\frac{\epsilon}{1+\epsilon}}{p r_0}+\frac{\frac{1}{1+\epsilon}}{pr_1}=\frac{1}{p} \q \text{ and }\q \frac{\frac{\epsilon}{1+\epsilon}}{p_j r_0}+\frac{\frac{1}{1+\epsilon}}{p_jr_1}=\frac{1}{p_j}~\text{ for }~j=1,\dots,m.$$
In order to use an interpolation in Lemma \ref{weightinterpol} with $\theta=\frac{1}{r_1(1+\epsilon)}$, we  see that  there exists $\delta_0>0$ such that
\begin{equation}\label{con1interest}
\big\Vert T_{K_{\mu}}(f_1,\dots,f_m)\big\Vert_{L^{p r_0}(\bbrn)}\lesssim 2^{-\delta_{0} \mu}\Vert \Omega\Vert_{L^q(\mathbb{S}^{mn-1})}\prod_{j=1}^{m}\Vert f_j\Vert_{L^{p_j r_0}(\bbrn)}
\end{equation}
for all $f_1,\dots,f_m\in\mathscr{C}_c^{\infty}(\bbrn)$,
which follows from \eqref{nonweightkeyest}.
Moreover, Lemma \ref{mweightincl} implies $\www^{1+\epsilon}\in \mathrm{A}_{r_1\ppp /q'}$ and thus Corollary \ref{qconjulessp} yields
\begin{equation}\label{con2interest}
\big\Vert T_{K_{\mu}}(g_1,\dots,g_m     )\big\Vert_{L^{pr_1}(v_{\www}^{1+\epsilon})}\lesssim  \mu \Vert \Omega\Vert_{L^q(\mathbb{S}^{mn-1})}\prod_{j=1}^{m}\Vert g_j\Vert_{L^{p_j r_1}(w_j^{1+\epsilon})}
\end{equation}
for all $g_j\in L^{\infty}_c(\bbrn)$.
Here, we note that
$p_1r_1,\dots,p_m r_1>q'$.

For any positive integers $M,N>0$  let
$$\mathcal{D}_M:=\big\{x\in \bbrn: v_{\www}(x)\le M \big\}$$
and
\begin{equation*}
\mathcal{E}_N^j:=\big\{x\in \bbrn: 1/N<w_j(x)<N \big\}, \q j=1,\dots,m.
\end{equation*}
For each $z\in \overline{\mathbf{S}}$ we define
\begin{align*}
V_z^{N}\big(f_1,\dots,f_m \big)(x):=v_{\www}(x)^{\frac{1+\epsilon}{pr_1}z}T_{K_{\mu}}\Big( \chi_{\mathcal{E}_N^1}f_1 w_1^{-\frac{1+\epsilon}{p_1r_1}z},\dots, \chi_{\mathcal{E}_N^m}f_m w_m^{-\frac{1+\epsilon}{p_mr_1}z}           \Big)(x)
\end{align*}
and
\begin{align*}
U_z^{M,N}\big(f_1,\dots,f_m \big)(x)&:=\chi_{\mathcal{D}_M}(x)V_z^{N}\big(f_1,\dots,f_m \big)(x)
\end{align*}
for all $f_1,\dots,f_m\in\mathscr{C}_c^{\infty}(\bbrn)$.
Then the estimates \eqref{con1interest} and  \eqref{con2interest}, with $g_j=\chi_{\mathcal{E}_N^j}f_j w_j^{-\frac{1+\epsilon}{p_jr_1}z}\in L^{\infty}_c(\bbrn)$, yield
\begin{align}
\big\Vert U_{iy}^{M,N}(f_1,\dots,f_m)\big\Vert_{L^{pr_0}(\bbrn)}&\lesssim 2^{-\delta_0 \mu} \Vert \Omega\Vert_{L^q(\mathbb{S}^{mn-1})} \prod_{j=1}^{m}\Vert f_j\Vert_{L^{p_jr_0}(\bbrn)}\label{preestinter1}\\
\big\Vert U_{1+iy}^{M,N}(f_1,\dots,f_m)\big\Vert_{L^{pr_1}(\bbrn)}&\lesssim \mu \Vert \Omega\Vert_{L^q(\mathbb{S}^{mn-1})} \prod_{j=1}^{m}\Vert f_j\Vert_{L^{p_jr_1}(\bbrn)},\label{preestinter2}
\end{align}
uniformly in $M,N,$ and $y$, respectively, for all $f_1,\dots,f_m\in\mathscr{C}_c^{\infty}(\bbrn)$.

In addition, $U_z^{M,N}$ is a multilinear operator mapping $\mathscr{C}_c^{\infty}(\bbrn)\times \cdots\times \mathscr{C}_c^{\infty}(\bbrn)$ to $L_{\mathrm{loc}}^1(\bbrn)$ and $\{U_z^{M,N}\}$ becomes an analytic family of admissible growth. 
Indeed, if $p\ge 1$, then H\"older's inequality and  Corollary \ref{qconjulessp}  prove that for any compact set $B$ in $\bbrn$,
\begin{align*}
&\int_B \big|  U_{z}^{M,N}(f_1,\dots,f_m)(x) \big|\; dx\\
&=\int_{B\cap \mathcal{D}_M}v_{\www}(x)^{\frac{1+\epsilon}{pr_1}\mathrm{Re(z)}}\Big| T_{K_{\mu}}\Big( \chi_{\mathcal{E}_N^1}f_1 w_1^{-\frac{1+\epsilon}{p_1r_1}z},\dots, \chi_{\mathcal{E}_N^m}f_m w_m^{-\frac{1+\epsilon}{p_mr_1}z}           \Big)(x) \Big|\; dx\\
&\le |B|^{1-\frac{1}{pr_1}}\Big\Vert   T_{K_{\mu}}\Big( \chi_{\mathcal{E}_N^1}f_1 w_1^{-\frac{1+\epsilon}{p_1r_1}z},\dots, \chi_{\mathcal{E}_N^m}f_m w_m^{-\frac{1+\epsilon}{p_mr_1}z}           \Big)  \Big\Vert_{L^{pr_1}(v_{\www}^{(1+\epsilon)\mathrm{Re}(z)})}\\
&\lesssim_B\mu \Vert \Omega\Vert_{L^q(\mathbb{S}^{mn-1})}\prod_{j=1}^{m}\Vert f_j\Vert_{L^{p_jr_1}(\bbrn)}
\end{align*}
uniformly in $\mathrm{Im}(z)$, $M$, and $N$, as $pr_1>1$, $p_1r_1,\dots,p_mr_1>q'$, and $\www^{(1+\epsilon)\mathrm{Re}(z)}\in \mathrm{A}_{r_1\ppp/q'}$ due to Lemma \ref{weightlowpore}.
If $p<1(<q')$, then $q'<p_1/p,\dots,p_m/p$, and
Lemma \ref{mweightincl} deduces $\www^{1+\epsilon}\in \mathrm{A}_{\ppp/(pq')}$. 
When $\frac{\mathrm{Re}(z)}{pr_1}>1$,
 Corollary \ref{qconjulessp} proves that for any compact set $B$ in $\bbrn$, 
\begin{align*}
&\int_B \big|  U_{z}^{M,N}(f_1,\dots,f_m)(x) \big|\; dx\\
&=\int_{B\cap \mathcal{D}_M}v_{\www}(x)^{\frac{1+\epsilon}{pr_1}\mathrm{Re(z)}}\Big| T_{K_{\mu}}\Big( \chi_{\mathcal{E}_N^1}f_1 w_1^{-\frac{1+\epsilon}{p_1r_1}z},\dots, \chi_{\mathcal{E}_N^m}f_m w_m^{-\frac{1+\epsilon}{p_mr_1}z}           \Big)(x) \Big|\; dx\\
&\lesssim_M \int_{\bbrn}\Big| T_{K_{\mu}}\Big( \chi_{\mathcal{E}_N^1}f_1 w_1^{-\frac{1+\epsilon}{p_1r_1}z},\dots, \chi_{\mathcal{E}_N^m}f_m w_m^{-\frac{1+\epsilon}{p_mr_1}z}           \Big)(x) \Big| v_{\www}(x)^{1+\epsilon} \; dx\\
&\lesssim \mu \Vert \Omega\Vert_{L^q(\mathbb{S}^{mn-1})}\prod_{j=1}^{m}\Big\Vert    \chi_{\mathcal{E}_N^j}f_j w_j^{-\frac{1+\epsilon}{p_jr_1}\mathrm{Re}(z)}    \Big\Vert_{L^{p_j/p}(w_j^{1+\epsilon})}\\
&\lesssim_N \mu \Vert \Omega\Vert_{L^q(\mathbb{S}^{mn-1})}\prod_{j=1}^{m}\Vert f_j\Vert_{L^{p_j/p}(\bbrn)}
\end{align*}
uniformly in $\mathrm{Im}(z)$.
When $\frac{\mathrm{Re}(z)}{pr_1}\le 1 $, Lemma \ref{weightlowpore} implies $\www^{\frac{1+\epsilon}{pr_1}\mathrm{Re}(z)}\in \mathrm{A}_{\ppp/(pq')}$. Then Corollary \ref{qconjulessp} proves that
\begin{align*}
&\int_{\bbrn} \big|  U_{z}^{M,N}(f_1,\dots,f_m)(x) \big|\; dx\\
&=\bigg\Vert     T_{K_{\mu}}\Big( \chi_{\mathcal{E}_N^1}f_1 w_1^{-\frac{1+\epsilon}{p_1r_1}z},\dots, \chi_{\mathcal{E}_N^m}f_m w_m^{-\frac{1+\epsilon}{p_mr_1}z}           \Big)     \bigg\Vert_{L^1(v_{\www}^{\frac{1+\epsilon}{pr_1}\mathrm{Re(z)}} )}\\
&\lesssim \mu \Vert \Omega\Vert_{L^q(\mathbb{S}^{mn-1})}\prod_{j=1}^{m}\Big\Vert    \chi_{\mathcal{E}_N^j}f_j w_j^{-\frac{1+\epsilon}{p_jr_1}\mathrm{Re}(z)}    \Big\Vert_{L^{p_j/p}(w_j^{\frac{1+\epsilon}{pr_1}\mathrm{Re}(z)})}\\
&\lesssim \mu \Vert \Omega\Vert_{L^q(\mathbb{S}^{mn-1})}\prod_{j=1}^{m}\Vert f_j\Vert_{L^{p_j/p}(\bbrn)}
\end{align*}
uniformly in $\mathrm{Im}(z)$, $M$, and $N$.
Therefore, $U_{z}^{M,N}$ satisfies the assumptions on $T_z$ in Lemma \ref{weightinterpol}.
 Now applying the interpolation method in Lemma \ref{weightinterpol} between \eqref{preestinter1} and \eqref{preestinter2}, we obtain
$$\Big\Vert U_{\frac{1}{r_1(1+\epsilon)}}^{M,N}(f_1,\dots,f_m)\Big\Vert_{L^{p}(\bbrn)}\lesssim 2^{-\frac{\delta_0\epsilon}{1+\epsilon} \mu}\mu^{\frac{1}{1+\epsilon}} \Vert \Omega\Vert_{L^q(\mathbb{S}^{mn-1})} \prod_{j=1}^{m}\Vert f_j \Vert_{L^{p_j}(\bbrn)}$$
for $f_1,\dots,f_m\in\mathscr{C}_c^{\infty}(\bbrn)$,
where the constant in the inequality is independent of $M$ and $N$.
Since
$$\lim_{M\to \infty}U_{\frac{1}{r_1(1+\epsilon)}}^{M,N}\big(f_1,\dots,f_m\big)(x)         = V_{\frac{1}{r_1(1+\epsilon)}}^{N}\big(f_1,\dots,f_m\big)(x),$$
Fatou's lemma yields that
\begin{align}\label{tkmuchiejnest}
\Big\Vert T_{K_{\mu}}\Big(\chi_{\mathcal{E}_N^1}f_1 w_1^{-\frac{1}{p_1}},\dots, \chi_{\mathcal{E}_N^m}f_m w_m^{-\frac{1}{p_m}}    \Big) \Big\Vert_{L^p(v_{\www})}&=\Big\Vert V_{\frac{1}{r_1(1+\epsilon)}}^{N}(f_1,\dots,f_m)\Big\Vert_{L^{p}(\bbrn)}\nonumber\\
&\lesssim 2^{-\frac{\delta_0\epsilon}{1+\epsilon} \mu}\mu^{\frac{1}{1+\epsilon}} \Vert \Omega\Vert_{L^q(\mathbb{S}^{mn-1})} \prod_{j=1}^{m}\Vert f_j \Vert_{L^{p_j}(\bbrn)}
\end{align}
uniformly in $N$, for all $f_1,\dots,f_m\in\mathscr{C}_c^{\infty}(\bbrn)$.
Since $L^{\infty}_c(\bbrn)\subset L^{p_j}(\bbrn)$ and $\mathscr{C}_c^{\infty}(\bbrn)$ is dense in $L^{p_j}(\bbrn)$,
the estimate \eqref{tkmuchiejnest} still holds for $g_1,\dots,g_m\in L_c^{\infty}(\bbrn)$.
For any $f_j\in\mathscr{C}_c^{\infty}(\bbrn)$, setting
$$g_j= f_j w_j^{\frac{1}{p_j}}\chi_{\mathcal{E}_N^j}\in L_c^{\infty}(\bbrn),$$
we have
\begin{align*}
\Big\Vert T_{K_{\mu}}\big(\chi_{\mathcal{E}_N^1}f_1,\dots, \chi_{\mathcal{E}_N^m}f_m \big) \Big\Vert_{L^p(v_{\www})}&\lesssim 2^{-\frac{\delta_0\epsilon}{1+\epsilon} \mu}\mu^{\frac{1}{1+\epsilon}} \Vert \Omega\Vert_{L^q(\mathbb{S}^{mn-1})} \prod_{j=1}^{m}\Vert f_j \Vert_{L^{p_j}(w_j)}
\end{align*}
where the constant in the inequality is independent of $N$.
Now we claim that for each $f_j \in \mathscr{C}_c^{\infty}(\bbrn)$, $1\le j \le m$,
\begin{equation}\label{lpvwconverge}
T_{K_{\mu}}\big(\chi_{\mathcal{E}_N^1}f_1 ,\dots, \chi_{\mathcal{E}_N^m}f_m \big)  \to T_{K_{\mu}}\big(f_1 ,\dots, f_m \big)  ~\text{ in } ~L^p(v_{\www})
\end{equation}
as $N\to \infty$. Indeed, applying the estimate \eqref{tkmug1gm2mumn2d},
\begin{align*}
&\bigg\Vert T_{K_{\mu}}\big(\chi_{\mathcal{E}_N^1}f_1 ,\dots, \chi_{\mathcal{E}_N^m}f_m \big)  - T_{K_{\mu}}\big(f_1 ,\dots, f_m \big) \bigg\Vert_{L^p(v_{\www})}\\
&\le \Big\Vert T_{K_{\mu}}\big( (\chi_{\mathcal{E}_N^1}-1)f_1,\chi_{\mathcal{E}_N^2}f_2,\dots, \chi_{\mathcal{E}_N^m}f_m      \big)\Big\Vert_{L^p(v_{\www})}\\
&\q+  \Big\Vert T_{K_{\mu}}\big(f_1, (\chi_{\mathcal{E}_N^2}-1)f_2,\chi_{\mathcal{E}_N^3}f_3,\dots, \chi_{\mathcal{E}_N^m}f_m     \big)\Big\Vert_{L^p(v_{\www})}\\
&\qq\qq\qq\qq\qq\qq\vdots\\
&\q +\Big\Vert T_{K_{\mu}}\big(f_1,\dots, f_{m-1},    (\chi_{\mathcal{E}_N^m}-1)f_m   \big)\Big\Vert_{L^p(v_{\www})}\\
&\lesssim 2^{\mu(mn+2-\delta)} \Vert \Omega\Vert_{L^q(\mathbb{S}^{mn-1})}\sum_{k=1}^{m}\Big\Vert \big( \chi_{\mathcal{E}_N^k}-1\big)f_k\Big\Vert_{L^{p_k}(w_k)}\prod_{j\not= k}\Vert f_j\Vert_{L^{p_j}(w_j)}.
\end{align*}
Here, we may also employ Corollary \ref{qconjulessp}, instead of \eqref{tkmug1gm2mumn2d}, in the case when $p_1,\dots,p_m>q'$.
The right-hand side converges to zero as $N\to \infty$ because the dominated convergence theorem yields that
\begin{align*}
\lim_{N\to \infty}\big\Vert \big( \chi_{\mathcal{E}_N^k}-1\big)f_k\big\Vert_{L^{p_k}(w_k)}= \big\Vert f_k \lim_{N\to \infty}\big( \chi_{\mathcal{E}_N^k}-1\big)\big\Vert_{L^{p_k}(w_k)}=0
\end{align*}
as $N\to \infty$. This proves \eqref{lpvwconverge}.
Hence, we have
\begin{align*}
\big\Vert T_{K_{\mu}}\big(f_1 ,\dots, f_m    \big) \big\Vert_{L^p(v_{\www})}\lesssim 2^{-\frac{\delta_0\epsilon}{1+\epsilon} \mu}\mu^{\frac{1}{1+\epsilon}} \Vert \Omega\Vert_{L^q(\mathbb{S}^{mn-1})} \prod_{j=1}^{m}\Vert f_j \Vert_{L^{p_j}(w_j)}
\end{align*}
for all $f_j \in \mathscr{C}_c^{\infty}(\bbrn)$, $1\le j\le m$,
and then by taking the sum over $\mu>0$, the desired result follows.

\hfill

\section{Proof of Proposition \ref{mainpropo}}\label{mainpropoprsec}
Let $g_1,\dots,g_m\in L^{\infty}_c(\bbrn)$.
We need to prove that for each fixed cube $Q$ containing $x$,
\begin{align}\label{pro9mainest}
&\inf_{c_Q\in\mathbb{C}}\bigg( \frac{1}{|Q|}\int_{Q} \big| T_{K_{\mu}}\big(g_1,\dots,g_m\big)(y)-c_Q\big|^r \; dy\bigg)^{1/r}\nonumber\\
& \lesssim  \mu\Vert \Omega\Vert_{L^q(\mathbb{S}^{mn-1})} \mathbf{M}_{q'}\big(g_1,\dots,g_m\big)(x)
\end{align}
uniformly in $Q$.

Let $Q^*$ be the concentric dilate of $Q$ with $\ell(Q^*)= 10\sqrt{mn}\ell(Q)$ and let $\chi_{Q^*}$ indicate the characteristic function of $Q^*$. We divide
$$g_j=g_j\chi_{Q^*}+g_j\chi_{(Q^*)^c}=:g_j^{(0)}+g_j^{(1)} \q \text{ for }~ j=1,\dots,m$$
so that $T_{K_{\mu}}(g_1,\dots,g_m)$ can be expressed as the sum of $2^m$ different terms of the form
$$T_{K_{\mu}}\big(g_1^{(\lambda_1)},g_2^{(\la_2)}\dots,g_m^{(\la_m)}\big)$$
where each $\lambda_j$ is $0$ or $1$. Here, we note that each $g_j^{(\lambda_j)}$ also belongs to $L^{\infty}_c(\bbrn)$ so that the last expression is well-defined almost everywhere, in view of \eqref{tkmug1gmptest}.
Then the left-hand side of \eqref{pro9mainest} is bounded by the sum of 
$$\mathcal{I}_1:=\bigg(\frac{1}{|Q|}\int_Q \big| T_{K_{\mu}}\big(g_1^{(0)},\dots,g_m^{(0)}\big)(y)\big|^r\; dy \bigg)^{1/r}$$
and
$$\mathcal{I}_2:=\inf_{c_Q\in\bbc}\bigg(\frac{1}{|Q|}\int_Q \Big|\Big( T_{K_{\mu}}\big(g_1,\dots,g_m\big)(y)-T_{K_{\mu}}\big(g_1^{(0)},\dots,g_m^{(0)}\big)(y)    \Big)   -c_Q\Big|^r\; dy\bigg)^{1/r}.$$

Using \eqref{bcrnlpest} in the setting $p_1=\cdots=p_m=mr>1$, which satisfies \eqref{qcondequi},
\begin{align*}
\mathcal{I}_1&\le \frac{1}{|Q|^{1/r}}\big\Vert T_{K_{\mu}}\big( g_1^{(0)},\dots,g_m^{(0)}\big) \big\Vert_{L^r(\bbrn)}\\
&\lesssim 2^{-\delta_0\mu}\Vert \Omega\Vert_{L^q(\mathbb{S}^{mn-1})}\frac{1}{|Q|^{1/r}}\prod_{j=1}^{m}\big\Vert g_j\big\Vert_{L^{mr}(Q^*)}\\
&\lesssim 2^{-\delta_0\mu}\Vert \Omega\Vert_{L^q(\mathbb{S}^{mn-1})}\mathbf{M}_{q'}\big(g_1,\dots,g_m\big)(x).
\end{align*}

In order to estimate $\mathcal{I}_2$, we employ the decomposition \eqref{kgmdef} and take
$$c_Q=\sum_{\gamma\in\bbz:2^{\gamma}\ell(Q)<2^{-(mn+1)\mu}}\Big( T_{K_{\mu}^{\gamma}}\big( g_1,\dots,g_m\big)(x)-  T_{K_{\mu}^{\gamma}}\big( g_1^{(0)},\dots,g_m^{(0)}\big)(x)    \Big)$$
to deduce
\begin{align*}
&\Big|\Big( T_{K_{\mu}}\big(g_1,\dots,g_m\big)(y)-T_{K_{\mu}}\big(g_1^{(0)},\dots,g_m^{(0)}\big)(y)    \Big)   -c_Q\Big|\\
&\le 
\sum_{\substack{\lambda_1,\dots,\lambda_m\in\{0,1\}\\ (\lambda_1,\dots,\lambda_m)\not= (0,\dots,0)}} \bigg( \sum_{\gamma\in\bbz:2^{\gamma}\ell(Q)\ge 1}\big| T_{K_{\mu}^{\gamma}}\big( g_1^{(\lambda_1)},\dots,g_m^{(\lambda_m)}\big)(y)\big| \\
&\q\q   +\sum_{\gamma\in\bbz:       2^{-\mu (mn+1)}\le 2^{\gamma}\ell(Q)<1}     \big| T_{K_{\mu}^{\gamma}}\big( g_1^{(\lambda_1)},\dots,g_m^{(\lambda_m)}\big)(y)\big|      \\
&\q\q\q +   \sum_{\gamma\in\bbz:2^{\gamma}\ell(Q)<2^{-\mu (mn+1)}}\Big|  T_{K_{\mu}^{\gamma}}\big( g_1^{(\lambda_1)},\dots,g_m^{(\lambda_m)}\big)(y)-  T_{K_{\mu}^{\gamma}}\big(  g_1^{(\lambda_1)},\dots,g_m^{(\lambda_m)}\big)(x) \Big| \bigg)\\
&=:\sum_{\substack{\lambda_1,\dots,\lambda_m\in\{0,1\}\\ (\lambda_1,\dots,\lambda_m)\not= (0,\dots,0)}} \Big(\Gamma_{Q,1}^{(\lambda_1,\dots,\lambda_m)}(y) +\Gamma_{Q,2}^{(\lambda_1,\dots,\lambda_m)}(y) +\Gamma_{Q,3}^{(\lambda_1,\dots,\lambda_m)}(y) \Big).
\end{align*}
Now we claim that for each $(\lambda_1,\dots,\lambda_m)\not= (0,\dots,0)$,
\begin{equation}\label{gq1l1myest}
\bigg( \frac{1}{|Q|}\int_Q  \Big| \Gamma_{Q,1}^{(\la_1,\dots,\la_m)} (y) \Big|^r   \; dy  \bigg)^{1/r}\lesssim \Vert \Omega\Vert_{L^1(\mathbb{S}^{mn-1})} \mathbf{M}\big( g_1,\dots,g_m\big)(x),
\end{equation}
\begin{equation}\label{gq2l1myest}
\bigg( \frac{1}{|Q|}\int_Q  \Big| \Gamma_{Q,2}^{(\la_1,\dots,\la_m)} (y) \Big|^r   \; dy  \bigg)^{1/r}\lesssim \mu \Vert \Omega\Vert_{L^q(\mathbb{S}^{mn-1})}\mathbf{M}_{q'}\big(g_1,\dots,g_m\big)(x),
\end{equation}
\begin{equation}\label{gq3l1myest}
\bigg( \frac{1}{|Q|}\int_Q  \Big| \Gamma_{Q,3}^{(\la_1,\dots,\la_m)} (y) \Big|^r   \; dy  \bigg)^{1/r}\lesssim \Vert \Omega\Vert_{L^1(\mathbb{S}^{mn-1})} \mathbf{M}\big( g_1,\dots,g_m\big)(x)
\end{equation}
uniformly in $Q$.
Then \eqref{pro9mainest} obviously follows.

For the proof of \eqref{gq1l1myest}, \eqref{gq2l1myest}, and \eqref{gq3l1myest}, permuting the variables, without loss of generality, we may assume $\la_1=\cdots=\la_{\kappa}=0$, $\la_{\kappa+1}=\cdots=\la_m=1$ for $1\le \kappa\le m-1$.
For simplicity, we use the notation
$$\mathcal{P}_Q^{\kappa}:=\big\{ \uuu\in (\bbrn)^m: u_1,\dots,u_{\kappa}\in Q^*,\, u_{\kappa+1},\dots,u_m\in (Q^*)^c\big\}.$$

\hfill

\subsection{Proof of \eqref{gq1l1myest}}

If $x,y\in Q$, $2^{-\gamma-1}\le |\zzz|\le 2^{-\gamma+1}$, $\uuu\in \mathcal{P}_Q^{\kappa}$, and $\ell(Q)\ge 2^{-\gamma}$,
then
\begin{align*}
&\big|\big(y-z_1-u_1,y-z_2-u_2,\dots,y-z_m-u_m \big) \big|\\
&\ge \big|\big(0,\dots,0,y-z_{\kappa+1}-u_{\kappa+1},\dots,y-z_m-u_m \big) \big|\\
&\ge \big|\big(0,\dots,0,x-u_{\kappa+1},\dots,x-u_m \big) \big|\\
&\qq\qq - \big| \big(0,\dots,0,x-y,\dots,x-y\big)\big|-|(0,\dots,0,z_{\kappa+1},\dots,z_m)|\\
&\gtrsim \big|\big(0,\dots,0,x-u_{\kappa+1},\dots,x-u_m \big) \big|\\
&\gtrsim \big|\big(x-u_1,\dots,x-u_m \big) \big|\gtrsim \ell(Q).
\end{align*}
Therefore, for $x,y\in Q$ and $2^{-\gamma-1}\le |\zzz|\le 2^{-\gamma+1}$,
\begin{align*}
&\big| \Psi_{\mu+\gamma}\ast \big( g_1^{(0)}\otimes\cdots\otimes g_{\kappa}^{(0)}\otimes g_{\kappa+1}^{(1)}\otimes\cdots\otimes g_{m}^{(1)}\big)(y-z_1,\dots,y-z_m)\big|\\
&\lesssim_L \frac{1}{2^{\mu+\gamma}\ell(Q)}\int_{(\bbrn)^m}\frac{2^{(\mu+\gamma)mn}}{(1+2^{\mu+\gamma}|(x-u_1,\dots,x-u_m)|)^L}\prod_{j=1}^{m}\big| g_j(u_j)\big|\; d\uuu\\
&\lesssim \frac{1}{2^{\mu+\gamma}\ell(Q)}\mathbf{M}\big( g_1,\dots,g_m\big)(x)
\end{align*} 
for $L>mn$, where the last inequality follows from Lemma \ref{multmaximalpt}.
 This yields that for $x,y\in Q$,
\begin{align*}
&\Big| T_{K_{\mu}^{\gamma}}\big(g_1^{(0)},\dots,g_{\kappa}^{(0)},g_{\kappa+1}^{(1)},\dots,g_{m}^{(1)}\big)(y)\Big|\\
&\le \int_{2^{-\gamma-1}\le |\zzz|\le 2^{-\gamma+1}}\big| K^{\gamma}(\zzz)\big|\big| \Psi_{\mu+\gamma}\ast \big( g_1^{(0)}\otimes\cdots\otimes g_{\kappa}^{(0)}\otimes g_{\kappa+1}^{(1)}\otimes\cdots\otimes g_{m}^{(1)}\big)(y-z_1,\dots,y-z_m)\big|\; d\zzz\\
&\lesssim  \frac{1}{2^{\mu+\gamma}\ell(Q)}\big\Vert K^{\gamma}\big\Vert_{L^1((\bbrn)^m)} \mathbf{M}\big( g_1,\dots,g_m\big)(x)\\
&\lesssim  \frac{1}{2^{\mu+\gamma}\ell(Q)}\Vert \Omega\Vert_{L^1(\mathbb{S}^{mn-1})} \mathbf{M}\big( g_1,\dots,g_m\big)(x).
\end{align*}
Finally, if $\mu>0$, then
\begin{align*}
\bigg( \frac{1}{|Q|}\int_Q  \Big| \Gamma_{Q,1}^{(\la_1,\dots,\la_m)} (y) \Big|^r   \; dy  \bigg)^{1/r}&\lesssim  \Vert \Omega\Vert_{L^1(\mathbb{S}^{mn-1})} \mathbf{M}\big( g_1,\dots,g_m\big)(x),
\end{align*}
as desired.

\hfill

\subsection{Proof of \eqref{gq2l1myest}}
We first claim that if $\ell(Q)< 2^{-\gamma}$ and $x,y\in Q$, then
\begin{equation}\label{tkmgf1fmy}
\Big| T_{K_{\mu}^{\gamma}}\big( g_1^{(\lambda_1)},\dots,g_m^{(\lambda_m)}\big)(y)\Big|\lesssim \Vert \Omega\Vert_{L^q(\mathbb{S}^{mn-1})}\mathbf{M}_{q'}\big(g_1,\dots,g_m\big)(x).
\end{equation}
Indeed, the left-hand side of \eqref{tkmgf1fmy} can be bounded by the sum of 
\begin{align*}
\mathcal{J}_1(y)&:=  \int_{|(x-u_1,\dots,x-u_m)|\le 10mn 2^{-\gamma}}\big|K_{\mu}^{\gamma}(y-u_1,\dots,y-u_m) \big|\prod_{j=1}^{m}|g_j(u_j)|\; d\uuu, \\
 \mathcal{J}_2(y)&:= \int_{|(x-u_1,\dots,x-u_m)|> 10mn 2^{-\gamma}}\big|K_{\mu}^{\gamma}(y-u_1,\dots,y-u_m) \big|\prod_{j=1}^{m}|g_j(u_j)|\; d\uuu. 
\end{align*}
By H\"older's inequality,
\begin{align*}
\mathcal{J}_1(y)&\le \big\Vert K_{\mu}^{\gamma}\big\Vert_{L^q((\bbrn)^m)} \bigg( \int_{|(x-u_1,\dots,x-u_m)|\le 10mn 2^{-\gamma}}\prod_{j=1}^{m}|g_j(u_j)|^{q'}\; d\uuu \bigg)^{1/q'}\\
&\lesssim 2^{-\gamma mn/q'}\big\Vert K^{\gamma}\big\Vert_{L^q((\bbrn)^m)}\mathbf{M}_{q'}\big(g_1,\dots,g_m\big)(x)\\
&\lesssim \Vert \Omega\Vert_{L^q(\mathbb{S}^{mn-1})}\mathbf{M}_{q'}\big(g_1,\dots,g_m\big)(x).
\end{align*}
Moreover, if $|(x-u_1,\dots,x-u_m)|>10 mn 2^{-\gamma}$, $2^{-\gamma-1}\le |\zzz|\le 2^{-\gamma+1}$, and $\ell(Q)<2^{-\gamma}$, then
\begin{align*}
\big| (y-u_1-z_1,\dots,y-u_m-z_m)\big|&\ge \big|(x-u_1,\dots,x-u_m) \big|-\sqrt{m}|x-y|-|\zzz|\\
&\gtrsim \big|(x-u_1,\dots,x-u_m) \big|
\end{align*}
and thus
\begin{align*}
\big| \Psi_{\mu+\gamma}(y-u_1-z_1,\dots,y-u_m-z_m)\big|\lesssim_L \frac{2^{(\mu+\gamma)mn}}{(1+2^{\mu+\gamma}|(x-u_1,\dots,x-u_m)|)^L}
\end{align*}
for $L>mn$.
Then Lemma \ref{multmaximalpt} deduces that
\begin{align*}
\mathcal{J}_2(y)&\le  \int_{(\bbrn)^m}\big|K^{\gamma}(\zzz) \big|\int_{(\bbrn)^m} \frac{2^{(\mu+\gamma)mn}}{(1+2^{\mu+\gamma}|(x-u_1,\dots,x-u_m)|)^L}\prod_{j=1}^{m}|g_j(u_j)|\; d\uuu \,d\zzz\\
&\lesssim \Vert K^{\gamma}\Vert_{L^1((\bbrn)^m)}\mathbf{M}\big(g_1,\dots,g_m\big)(x)\\
&\lesssim \Vert \Omega\Vert_{L^1(\mathbb{S}^{mn-1})}\mathbf{M}\big(g_1,\dots,g_m\big)(x),
\end{align*}
which completes the proof of \eqref{tkmgf1fmy}.

Finally, the left-hand side of \eqref{gq2l1myest} is controlled by a constant times
\begin{align*}
&\sum_{\gamma\in\bbz:       2^{-\mu (mn+1)}\le 2^{\gamma}\ell(Q)<1} \Vert \Omega\Vert_{L^q(\mathbb{S}^{mn-1})}\mathbf{M}_{q'}\big(g_1,\dots,g_m\big)(x)\\
&\lesssim_{m,n} \mu  \Vert \Omega\Vert_{L^q(\mathbb{S}^{mn-1})}\mathbf{M}_{q'}\big(g_1,\dots,g_m\big)(x),
\end{align*}
as desired.

\hfill

\subsection{Proof of \eqref{gq3l1myest}}
We first write
\begin{align*}
&\Big| T_{K_{\mu}^{\gamma}}\big(g_1^{(0)},\dots,g_{\kappa}^{(0)},g_{\kappa+1}^{(1)},\dots,g_m^{(1)}\big)(y)-T_{K_{\mu}^{\gamma}}\big(g_1^{(0)},\dots,g_{\kappa}^{(0)},g_{\kappa+1}^{(1)},\dots,g_m^{(1)}\big)(x)\Big|\\
&\le \int_{(\bbrn)^m}\big|K^{\gamma}(\zzz) \big|\Big| \Psi_{\mu+\gamma}\ast \big(  g_1^{(0)}\otimes\cdots\otimes g_{\kappa}^{(0)}\otimes g_{\kappa+1}^{(1)}\otimes\cdots\otimes g_{m}^{(1)}       \big)(y-z_1,\dots,y-z_m)\\
&\qq\qq\qq-\Psi_{\mu+\gamma}\ast \big(  g_1^{(0)}\otimes\cdots\otimes g_{\kappa}^{(0)}\otimes g_{\kappa+1}^{(1)}\otimes\cdots\otimes g_{m}^{(1)}       \big)(x-z_1,\dots,x-z_m)\Big|\; d\zzz\\
&\le \int_{(\bbrn)^m}\big| K^{\gamma}(\zzz)\big|\bigg(\int_{\mathcal{P}_Q^{\kappa}}  \Big| \Psi_{\mu+\gamma}(y-z_1-u_1,\dots,y-z_m-u_m)\\
&\qq\qq\qq -\Psi_{\mu+\gamma}(x-z_1-u_1,\dots,x-z_m-u_m)\Big| \prod_{j=1}^{m}\big|g_j(u_j) \big|       \; d\uuu \bigg)         \; d\zzz\\
&=\int_{(\bbrn)^m}\big| K^{\gamma}(\zzz)\big|\bigg(\int_{\mathfrak{U}^1_{\gamma}} \cdots      \; d\uuu \bigg)         \; d\zzz+\int_{(\bbrn)^m}\big| K^{\gamma}(\zzz)\big|\bigg(\int_{\mathfrak{U}^2_{\gamma}} \cdots      \; d\uuu \bigg)         \; d\zzz\\
&=:\mathfrak{T}^1_{\mu,\gamma}(y)+\mathfrak{T}^2_{\mu,\gamma}(y)
\end{align*}
where
\begin{align*}
\mathfrak{U}_{\gamma}^{1}&:=\big\{\uuu\in \mathcal{P}_Q^{\kappa}: |x-u_j|< 10\sqrt{mn}2^{-\gamma}~\text{ for all }~j=\kappa+1,\dots,m\big\}\\
\mathfrak{U}_{\gamma}^2&:=\big\{\uuu\in \mathcal{P}_Q^{\kappa}: |x-u_j|\ge 10\sqrt{mn}2^{-\gamma}~\text{ for some }~j=\kappa+1,\dots,m\big\}
\end{align*}
so that $\mathcal{P}_Q^{\kappa}$ can be expressed as the disjoint union of $\mathfrak{U}_{\gamma}^{1}$ and $\mathfrak{U}_{\gamma}^{2}$.

To deal with $\mathfrak{T}_{\gamma}^{1}$, we see that for $x,y\in Q$
\begin{align*}
& \Big| \Psi_{\mu+\gamma}(y-z_1-u_1,\dots,y-z_m-u_m) -\Psi_{\mu+\gamma}(x-z_1-u_1,\dots,x-z_m-u_m)\Big|\\
&\lesssim 2^{\mu+\gamma}|x-y|2^{(\mu+\gamma)mn}\lesssim 2^{\mu+\gamma}\ell(Q)2^{(\mu+\gamma)mn}
\end{align*}
and thus 
\begin{align}\label{fraktg1est}
\mathfrak{T}_{\gamma}^{1}&\lesssim 2^{(\mu+\gamma) mn}2^{\mu+\gamma}\ell(Q)\big\Vert K^{\gamma}\big\Vert_{L^1((\bbrn)^m)}\int_{\mathfrak{U}_{\gamma}^{1}} \prod_{j=1}^{m}\big|  g_j(u_j)   \big|    \; d\uuu \nonumber\\
&\lesssim 2^{\mu mn}2^{\mu+\gamma}\ell(Q)\Vert \Omega\Vert_{L^1(\mathbb{S}^{mn-1})} 2^{\gamma mn}\int_{|\uuu|\lesssim 2^{-\gamma}}\prod_{j=1}^{m}\big|  g_j(x-u_j)   \big|  \; d\uuu\nonumber\\
&\lesssim 2^{\mu mn}2^{\mu+\gamma}\ell(Q)\Vert \Omega\Vert_{L^1(\mathbb{S}^{mn-1})}\mathbf{M}\big(g_1,\dots,g_m\big)(x)
\end{align}
where we observe that
$$\mathfrak{U}_{\gamma}^{1}\subset \big\{\uuu\in(\bbrn)^m: |\uuu|\le 10nm 2^{-\gamma} \big\}$$
as  $\mu>0$ and $\ell(Q)<2^{-\gamma}$.

Moreover, if $\uuu\in \mathfrak{U}_{\gamma}^{2}$, $2^{-\gamma-1}\le |\zzz|\le 2^{-\gamma+1}$, and $x,y\in Q$, then
\begin{align*}
& \Big| \Psi_{\mu+\gamma}(y-z_1-u_1,\dots,y-z_m-u_m) -\Psi_{\mu+\gamma}(x-z_1-u_1,\dots,x-z_m-u_m)\Big|\\
&\lesssim_L 2^{\mu+\gamma}|x-y| \int_0^1 \frac{2^{(\mu+\gamma)mn}}{(1+2^{\mu+\gamma}|(ty+(1-t)x-z_1-u_1,\dots,ty+(1-t)x-z_m-u_m)|)^L}\; dt\\
&\lesssim_L 2^{\mu+\gamma}\ell(Q)\frac{2^{(\mu+\gamma)mn}}{(1+2^{\mu+\gamma}|(x-u_1,\dots,x-u_m)|)^L}
\end{align*}
for $L>mn$,
since
\begin{align*}
&|(ty+(1-t)x-z_1-u_1,\dots,ty+(1-t)x-z_m-u_m)|\\
&\ge |(0,\dots,0,ty+(1-t)x-z_{\kappa+1}-u_{\kappa+1},\dots,ty+(1-t)x-z_m-u_m)|\\
&\ge |(0,\dots,0,x-u_{\kappa+1},\dots,x-u_m)|-\sqrt{m}|x-y|-|\zzz|\\
&\ge |(0,\dots,0,x-u_{\kappa+1},\dots,x-u_m)|-\sqrt{nm} \ell(Q)-2^{-\gamma+1}\\
&\ge  |(0,\dots,0,x-u_{\kappa+1},\dots,x-u_m)|-2\sqrt{nm}2^{-\gamma}\\
&\gtrsim   |(0,\dots,0,x-u_{\kappa+1},\dots,x-u_m)|\\
&\gtrsim  |(x-u_{1},\dots,x-u_m)|.
\end{align*}
This proves
\begin{align}\label{fraktg2est}
\mathfrak{T}_{\gamma}^{2}&\lesssim 2^{\mu+\gamma}\ell(Q)\big\Vert K^{\gamma}\big\Vert_{L^1((\bbrn)^m)} \int_{(\bbrn)^m}\frac{2^{(\mu+\gamma)mn}}{(1+2^{\mu+\gamma}|(x-u_1,\dots,x-u_m)|)^L}\prod_{j=1}^{m}\big| g_j(u_j)\big|\; d\uuu\nonumber\\
&\lesssim  2^{\mu+\gamma}\ell(Q)\Vert \Omega\Vert_{L^1(\mathbb{S}^{mn-1})}\mathbf{M}\big(g_1,\dots,g_m\big)(x).
\end{align}

Finally, \eqref{gq3l1myest} follows from \eqref{fraktg1est} and \eqref{fraktg2est}.

\hfill

\medskip
\noindent {\bf Acknowledgment:} Part of this work was carried out during my research stay at OIST. I would like to thank Xiaodan Zhou and Qing Liu for their invitation and hospitality during my stay.


\begin{thebibliography}{99}





\bibitem{Ca_Zy1952}
A. P. Calder\'on and A. Zygmund, \emph{On the existence of certain singular integrals},  Acta Math.  \textbf{88} (1952),  85--139.


\bibitem{Ca_Zy1956}
A. P. Calder\'on and A. Zygmund, \emph{On singular integrals},   
Amer. J. Math.  \textbf{78} (1956),  289--309.


\bibitem{Ch_Ru1988}
M. Christ and J.-L. Rubio de Francia ,  \emph{Weak type $(1,1)$ bounds for rough operators II}, Invent. Math. \textbf{93} (1988), 225--237.




\bibitem{Co_Me1975} 
R. R. Coifman and Y. Meyer,  
\emph{On commutators of singular integrals and bilinear singular integrals}, 
 Trans. Amer. Math. Soc. {\bf 212} (1975), 315--331.



\bibitem{Co_We1977}
R. R. Coifman and G. Weiss,  \emph{Extensions of Hardy spaces and their use in analysis}, Bull. Amer. Math. Soc. \textbf{83} (1977), 569-645.




\bibitem{Co_Cu_Di_Ou2017}
J. M. Conde-Alonso, A. Culiuc, F. Di Plinio, and Y. Ou,  \emph{A sparse domination principle for rough singular integrals}, Anal. PDE \textbf{10} (2017), 1255-1284.




\bibitem{Co1979}
W. C. Connett,  \emph{Singular integrals near $L^1$}, in Harmonic analysis in Euclidean spaces, Part 1 (Williamstown 1978), Proc. Sympos. Pure Math. \textbf{35} (1979), 163 --165.


\bibitem{Da1988}
I. Daubechies,  \emph{Orthonormal bases of compactly supported wavelets}, Comm. Pure Appl. Math. \textbf{41}
(1988), 909--996.


 \bibitem{Do_Sl2024} G. Dosidis and L. Slav\'ikov\'a, 
\emph{Multilinear singular integrals with homogeneous kernels near $L^1$}, 
 Math. Ann. \textbf{389} (2024), 2259--2271. 


\bibitem{Duan1993}
J. Duoandikoetxea,  \emph{Weighted norm inequalities for homogeneous singular integrals}, Trans. Amer. Math. Soc. \textbf{336} (1993), 869--880.



\bibitem{Du_Ru1986}
J. Duoandikoetxea and J.-L. Rubio de Francia,  \emph{Maximal and singular integral operators via Fourier transform estimates}, Invent. Math. \textbf{84} (1986), 541--561.




\bibitem{Fe_St1972}
C. Fefferman and E. M. Stein,  \emph{$H^p$ spaces of several variables}, Acta Math. \textbf{129}
(1972) 137--193.

\bibitem{Ga_Ru1985}
J. Garcia-Cuerva and J.-L. Rubio de Francia,  \emph{Weighted Norm Inequalities and Related Topics}, North-Holland Math. Studies, Vol. 116, North-Holland, Amsterdam, 1985.






\bibitem{Gr_He_Ho2018}
L. Grafakos, D. He, and P. Honz\'ik, \emph{Rough bilinear singular integrals}, Adv. Math. \textbf{326} (2018), 54--78.


\bibitem{Gr_He_Ho_Park2023}
L. Grafakos, D. He,   P. Honz\'ik, and B. Park, \emph{Initial $L^2\times\cdots\times L^2 $ bounds for   multilinear operators},  Trans. Amer. Math. Soc. \textbf{376} (2023), 3445-3472. 


\bibitem{Gr_He_Ho_Park_JLMS}
L. Grafakos, D. He, P. Honz\'ik, and B. Park,  \emph{Multilinear rough singular integral operators}, J. London Math. Soc. \textbf{109} (2024), e12867, 35pp.

\bibitem{Gr_He_Sl2020}
L. Grafakos, D. He, and L. Slav\'ikov\'a, \emph{$L^2\times L^2\to L^1$ boundedness criteria}, Math. Ann. \textbf{376} (2020), 431--455.



\bibitem{Gr_Ou2022}
L. Grafakos and E. M. Ouhabaz,  \emph{Interpolation for analytic families of multilinear operators on metric measure spaces}, Studia Math. \textbf{267} (2022), 37--57. 



\bibitem{GT-Indiana} 
L. Grafakos and R. H. Torres, 
\emph{Maximal operator and weighted norm inequalities for multilinear singular integrals}, 	Indiana Univ.  Math. J.  \textbf{51} (2002),   1261--1276.

\bibitem{Gr_To2002} 
L. Grafakos and R. H. Torres, 
\emph{Multilinear Calder\'on-Zygmund Theory}, Adv. Math.  \textbf{165} (2002),  124--164.


\bibitem{He_Park2023}
D. He and B. Park,  \emph{Improved estimates for bilinear rough singular integrals}, Math. Ann. \textbf{386} (2023), 1951-1978.



\bibitem{Ho1988}
S. Hofmann, \emph{Weak type $(1,1)$ boundedness of singular integrals with nonsmooth kernels}, Proc. Amer. Math. Soc. \textbf{103} (1988), 260--264.

\bibitem{Hy_Ro_Ta2017}
T. P. Hyt\"onen, L. Roncal, and O. Tapiola, \emph{Quantitative weighted estimates for rough homogeneous singular integrals}, Israel J. Math. \textbf{218} (2017), 133--164.



\bibitem{Le_Om_Pe_To_Tr2009}
A. Lerner, S. Ombrosi, C. P\'erez,  R. Torres, and R. Trujillo-Gonz\'alez,  \emph{New maximal functions and multiple weights for the multilinear Calder\'on-Zygmund theory}, Adv. Math. \textbf{220} (2009), 1222-1264.

\bibitem{Li_Pe_Ri_Ro2019}
K. Li, C. P\'erez, I.P. Rivera-R\'ios, and L. Roncal,  \emph{Weighted norm inequalities for rough singular integral operators}, J. Geom. Anal. \textbf{29} (2019), 2526-2564.



\bibitem{Li_Lo2001}
E.H. Lieb and M. Loss,  \emph{Analysis},  Graduate Studies in Mathematics, vol. 14, 2nd edition, American Mathematical Society, Providence 2001.


\bibitem{Mu1972} 
B. Muckenhoupt, \emph{Weighted norm inequalities for the Hardy maximal function}, Trans. Amer. Math. Soc. \textbf{165} (1972), 207--226.





\bibitem{Park_submitted2}
B. Park,  \emph{Multilinear estimates for maximal rough singular integrals}, Ann. Sc. Norm. Super. Pisa Cl. Sci., to appear.



\bibitem{Se1996}
A. Seeger, \emph{Singular integral operators with rough convolution kernels}, J. Amer. Math. Soc. \textbf{9} (1996), 95--105.

\bibitem{St2001}
A. Stefanov, \emph{Weak type estimates for certain Calder\'on-Zygmund singular integral operators}, Studia Math. \textbf{147} (2001), 1--13.


\bibitem{St1956}
E. M. Stein, \emph{Interpolation of linear operators}, Trans. Amer. Math. Soc. \textbf{83} (1956), 482--492.





\bibitem{Ta1999}
T. Tao, \emph{The weak type $(1,1)$ of $L\log{L}$ homogeneous convolution operators}, Indiana Univ. Math. J. \textbf{48} (1999), 1547--1548.

\bibitem{Va1996}
A. M. Vargas, \emph{Weighted weak type $(1,1)$ bounds for rough operators}, J. London Math. Soc. \textbf{54} (1996), 297--310.


\bibitem{Wa1990}
D. Watson, \emph{Weighted estimates for singular integrals via Fourier transform}, Duke Math. J. \textbf{60} (1990), 389--399.



\end{thebibliography}


\end{document}